\documentclass[10pt,reqno]{amsart}

\usepackage{amssymb}
\usepackage{amsfonts}
\usepackage{amsthm}
\usepackage{amsmath}
\usepackage{bbm}
\usepackage{xcolor}
\usepackage{enumerate}
\usepackage{mathabx}
\usepackage{lipsum,booktabs}
\usepackage{cite}
\usepackage[normalem]{ulem}
\usepackage{cleveref}

\DeclareRobustCommand{\rchi}{{\mathpalette\irchi\relax}}
\newcommand{\irchi}[2]{\raisebox{\depth}{$#1\chi$}} 

\numberwithin{equation}{section}
\allowdisplaybreaks

\newtheorem{theorem}{Theorem}[section]
\newtheorem{proposition}[theorem]{Proposition}

\newtheorem{lemma}[theorem]{Lemma}
\newtheorem{corollary}[theorem]{Corollary}

\theoremstyle{definition}
\newtheorem{definition}[theorem]{Definition}

\theoremstyle{remark}

\newcommand{\R}{\mathbb{R}}

\newcommand{\Z}{\mathbb{Z}}
\newcommand{\C}{\mathbb{C}}

\newcommand{\Rd}{\mathbb{R}^d}
\newcommand{\Rk}{\mathbb{R}^{d-1}}

\newcommand{\eps}{\varepsilon}

\newcommand{\scriptK}{\mathcal{K}}

\newcommand{\scriptM}{\mathcal{M}}

\newcommand{\scriptO}{\mathcal{O}}

\newcommand{\scriptX}{\mathcal{X}}

\newcommand{\qtq}[1]{\quad\text{#1}\quad}

\DeclareMathOperator*{\supp}{supp}

\begin{document}
\title{Compactness of a restricted X-ray transform}
\author{Chandan Biswas}

\thanks{This work was supported in part by NSF grants DMS-1266336 and DMS-1600458}
\address{Department of Mathematics, Indian Institute of Technology Bombay, Mumbai 400076}
\email{cbiswas@iitb.ac.in}

\begin{abstract}
We show that the X-ray transform with directions restricted along the moment curve possesses extremizers and that $L^p$-normalized extremizing sequences are precompact modulo symmetry. Our approach advances the Lorentz space method of Christ to a mixed norm Lebesgue space setting.
\end{abstract}
\date\today

\maketitle


\section{Introduction} 
Recently, the author proved that for any normalized sequence optimizing the norm of the averaging operator along the moment curve for an endpoint Lebesgue space bound, up to symmetry, either it converges in norm along a subsequence or after suitable scaling converges to an extremizer of the Xray transform with directions restricted along the moment curve. In this article, we establish the existence of extremizers and precompactness for this restricted Xray transform for all valid (including mixed norm) Lebesgue space estimates. To the best of our knowledge, this appears to be the first example of a precompactness result for an averaging operator in a mixed norm setting.

We introduce the operators. For a suitable, (say) smooth and compactly supported $f$ on $\Rd$, the averaging operator $T$ along the moment curve is defined by 
$$
Tf(x) := \int_\R f(x + (t, t^2, \ldots, t^d)) \, dt, \quad x \in \Rd
$$
and its X-ray transform restricted along the moment curve $X f$ is given by
\begin{equation}\label{E : def X}
Xf(t, y) := \int_{\R} f(s, y + s (t, t^2, \ldots, t^{d - 1})) \, ds, \quad (t, y) \in \R \times \Rk.
\end{equation}
Both $T, X$ extend as bounded operators from $L^{\frac{d^2 + d}{d^2 - d + 2}}$ into $L^{\frac{d + 1}{d - 1}}$ for $d \geq 3$ (\cite{Christ_1998, Littman1971, Oberlin1987, Stovall2009, ChristErdogan2002}).

Given a bounded operator $L : H \to K$ between two normed spaces $H, K$, we say a sequence $h_n \in H$ is extremizing if $\lim_n \frac{\|L h_n\|_K}{\|h_n\|_H} = \|L\|_{H \to K}$, and a nonzero $h \in H$ is an extremizer if $\|L h\|_K = \|L\|_{H \to K} \|h\|_H$. Additionally, by a symmetry of $L : H \to K$, we mean an isometry $S$ of $H$ for which there is an associated isometry $T_S$ of $K$ such that $L \circ S = T_S \circ L$. The following dichotomy was established in~\cite{CB}.

\begin{theorem}[~\cite{CB}]
Let $d \geq 3, p := \frac{d^2 + d}{d^2 - d + 2}$ and $q := \frac{d + 1}{d - 1}$. Let $\{f_n\}$ be normalized and extremizes $T : L^p \to L^q$. Then, either there exists a sequence of symmetries $\phi_n$ of $T$ such that, along a subsequence, $\phi_n(f_n)$ converges in norm to an extremizer of $T : L^p \to L^q$ or there exists a sequence of positive numbers $r_n$ so that the rescaled sequence $\{r_n^{d/p} f_n(r_n \cdot)\}$ converges in $L^p$ to an extremizer of $X : L^p \to L^q$.
\end{theorem}

In addition to the above, $X$ extends as a bounded operator (\cite{ChristErdogan2002, BS4, E, NL}) from $L^p$ into $L_t^q L_y^r$ for $(p, q, r) = (p_\theta, q_\theta, r_\theta)$ given by
\begin{equation}\label{E : def p_theta}
\Big(\frac{1}{p_\theta}, \frac{1}{q_\theta}, \frac{1}{r_\theta}\Big) := (1 -\theta) (1, 0, 1) + \theta (\frac d{d + 2}, \frac d{d + 2}, \frac{d^2 - d - 2}{d^2 + d - 2}), \quad \theta \in [0, 1),
\end{equation}
where the mixed norm $\|g\|_{L_t^q L_y^r}$ is given by
$$
\|g\|_{L_t^q L_y^r} := (\int_\R (\int_{\R^{d-1}} |g(t,y)|^r\,dy)^{q/r} dt)^{1/q}.
$$
We note that the non-mixed bound $X : L^{\frac{d^2 + d}{d^2 - d + 2}} \to L^{\frac{d + 1}{d - 1}}$ belongs to the family~\eqref{E : def p_theta} corresponding to
\begin{equation}\label{E : def theta_0}
\theta = \theta_0 := (d^2 + d - 2)/(d^2 + d), \qtq{so that} q_{\theta_0} = r_{\theta_0}.
\end{equation}

We now describe the symmetries of $X$. First, we note that for a diffeomorphism $\phi$ of $\Rd$, the map $\phi^* f := J_{\phi}^{\frac 1{p}} f \circ \phi$ is an isometry of $L^p(\Rd)$; likewise, for $\psi(t, y) = (\psi_1(t), \tilde \psi(y))$, a diffeomorphism of $\Rd$, the map $\psi^* g := |\psi_1'(t)|^{\frac 1q} J_{\tilde \psi}(y)^{\frac 1r} g \circ \psi$ is an isometry of $L_t^q L_y^r$. Three symmetries of $X$ in particular, translation: for $v \in \Rk$,
\begin{equation}\label{E : example symmetry}
(\phi_{(0, v)} (s, x), \psi_{(0, v)}(t, y)) := ((s, x + v), (t, y + v)),
\end{equation}
anisotropic scaling: with $S_\beta(x) := (\beta x_1, \beta^2 x_2, \ldots, \beta^{d - 1} x_{d - 1})$ and $\alpha \beta \neq 0$,
$$
(\phi_{(\alpha, \beta)}(s, x), \psi_{(\alpha, \beta)}(t, y)) := ((\alpha s, \alpha S_{\beta}(x)), (\beta t, \alpha S_{\beta}(y),
$$
and gliding along the moment curve $\gamma_0 (t) := (t, t^2, \ldots, t^{d - 1})$ : for fixed $s_0, t_0 \in \R$,
\begin{align*}
& \phi_{(s_0,t_0)}(s, x) := (s + s_0, G_{t_0}(x) + (s + s_0) \gamma_0(t_0))
\\
& \psi_{(s_0, t_0)}(t, y) := (t + t_0, G_{t_0}(y - s_0 \gamma_0(t))),
\end{align*}
and the group they generate are central to us. Here $G_s \in M_{d - 1}$ denotes the lower triangular matrix with $1$'s on the diagonal satisfying $\gamma_0 (s + t) = G_{s} \gamma_0 (t) + \gamma_0(s)$.

Our main result is the following.

\begin{theorem}\label{T : main thm}
Let $d \geq 3$ and $\theta \in [0, 1)$. Then there exist extremizers of $X : L^{p_\theta} \to L_t^{q_\theta} L_y^{r_\theta}$. Additionally, if $\theta \in (0, 1)$, for each normalized nonnegative extremizing sequence $\{f_n\}$ of $X : L^{p_\theta} \to L_t^{q_\theta} L_y^{r_\theta}$, there exists a sequence of symmetries $\{\phi_n^*\}$ such that along a subsequence $\{\phi_n^* f_n\}$ converges in norm to an extremizer.
\end{theorem}

\textbf{A bit of history of the operator $X$ :}
Given a nice, for example, compactly supported and smooth $f : \Rd \to \C$ and a line $l$ in $\Rd$, the value of the X-ray transform of $f$ at $l$ is defined to be the integral of $f$ over $l$. Introduced by Fritz~\cite{JF}, it is of paramount importance in integral geometry \cite{gel} and the study of PDE's (Indeed, it is known that the mixed norm Lebesgue space estimates of the X-ray transform implies the Kakeya Conjecture \cite{kattao, labatao, Wolf01}). Being overdetermined for $d\geq 3$, one naturally restricts the X-ray transform to the set of all lines whose directions are parametrized by a fixed curve $\gamma : \R \to \Rk$ and considers $X_\gamma$ which is defined by replacing the moment curve $\gamma_0$ in~\eqref{E : def X} by a fixed general curve $\gamma$.

Naturally, questions of regularity properties such as the Lebesgue space bounds for $X_\gamma$ have attracted continuous attention over several decades~\cite{ChristErdogan2008, GreenleafSeeger1998, E, Wolf01, Oberlin2000}. Despite this persistent effort, even for the model case $X = X_{\gamma_0}$, the question of identifying the range of exponents $(p, q, r)$ for which it maps $L^p$ into $L_t^q L_y^r$ is still open. Indeed, it has been conjectured~\cite[Conjecture $2.1$]{BS4} that~\eqref{E : def p_theta} gives the full range of exponents including at $\theta = 1$, and in~\cite{BS4} this has been verified away from the end point $\theta = 1$. Our main result, Theorem~\ref{T : main thm} sharpens these estimates.

The restricted Xray transforms are part of a larger class of singular integral operators known as the generalized radon transforms or averaging operators. These are defined by integrating a function over a family of lower dimensional submanifolds (see, e.g., ~\cite{Gressman2009, Tao-Wirght2001, NaiboJavier2013, Rubin2019, Rubinhyperbolic2019, Christ1984, Rubin2004, Rubin2023, Gressman2013, Gressmanintermediate2019, Gressman2007, BennettBezFlock2018, CarberySeegerWainger1999, CarberyChristWaingerWatson1989, CarberyRicciWright2003, PramanikSeeger2007, PramanikSeeger2006, PramanikSeeger2021, Gressman2025, Gressman2022, Stovall2011, ErdoganOberlin2010, ChristDendrinosStovallStreet2020} and the references therein). Sharp forms of these Lebesgue space bounds have been obtained in a few cases. For example, for the $k$-plane transform~\cite{Flockthesis, Flock2016, BernsteinLoss1997, BennettBezFlock2018, Rubin2019, Drouot2014, Druot2015}, averaging over hyperplanes~\cite{Christ2014}, averaging over the paraboloid~\cite{Christ_extremal, ChristXue2012}, Brascamp-Lieb inequalities~\cite{BennettCarberyChristTao2008}. These articles deal with non-mixed Lebesgue space estimates. One difference to note is that in this article the operator $X$ obeys mixed norm Lebesgue space bounds. Although we build on much of~\cite{Christ_quasiextremal, Christ_extremal, BS4} there are fundamental obstructions to adapt these arguments to our case. We give a brief account of this in the next section.

\section{An overview of the proof}
In Section~\ref{S : Proof of main thm}, we state two uniform localization results Propositions~\ref{P : freq loc} and~\ref{P : loc} for near extremizers. Assuming these, we give the proof of our main result Theorem~\ref{T : main thm}. This follows a standard proof in the sharp estimate literature.

We note a few preliminary results in Section~\ref{S : prelims}, which will be used repeatedly throughout this article. The extrapolation lemma~\cite[Lemma $6.3$]{BS4} is used to infer a bound on the range of $g$ for $(\rchi_E, g)$ a near extremizer pair (see~\eqref{E : def eps-quasi pair}), provided the function $g$ has constant $t$-fiber in~\eqref{E : def const $t$-fiber} which can be thought as an analogue of the product of a characteristic function in $t$ with a $G(y)$ in mixed norm spaces. The second result Proposition~\ref{P : similar vol E_1, E_2} is used to obtain uniform frequency estimates for near extremizers.

In Section~\ref{S : lorentz bound}, we obtain Lorentz space estimates for $X$. That Lorentz space estimates can be exploited to obtain uniform localizations for near extremizers was introduced in~\cite{Christ_extremal} and similar arguments have appeared in related contexts in non-mixed setting in~\cite{CB, DendrinosLaghiWright2009, NL, Stovall2010}. In the mixed norm setting, this argument has been adapted to prove a special case of bound~\eqref{E : eps-improv for theta large} when $f$ is a characteristic function (\cite[the first bound on page $25$]{BS4}). We extend this further and obtain Lorentz space estimates Theorems~\ref{T : L(p, s) Lorentz} and~\ref{T : L(q, s) lorentz} for $X$.

In Section~\ref{S : proof of fre loc}, we give the proof of the uniform frequency localization Proposition~\ref{P : freq loc}. At $\theta = \theta_0$ corresponding to the non-mixed bound $X : L^{p_{\theta_0}} \to L^{q_{\theta_0}}$ a straightforward adaptation of the argument in~\cite{Christ_extremal} suffices. However, at all other values of $\theta$, we face the fundamental obstruction that the mixed norm of a function does not necessarily decouple when decomposed into disjontly supported pieces, as evident from elementary examples, which plays a key role in the above non-mixed case (see bound~\eqref{E : small pure incidences}). To overcome this, given a pair $(E, F)$ and a positive $\lambda$, we introduce the notion of $\lambda$-strong points with large fiber (see~\eqref{E : def large fiber}), and combine this with the strict convexity argument in~\cite{Stovall2020} and our Lorentz space estimates. This is our primary innovation in this article. As this construction is of a general nature, we anticipate that this might be of interest for general mixed norm estimates.

We now move on to the proof of uniform localization for near extremizers Proposition~\ref{P : loc}. We begin this in Section~\ref{S : paraballs} by identifying a natural collection of quasiextremal pairs of $X$ (see~\eqref{E : def eps-quasi pair}) known as paraballs in literature and show that a general quasiextremal pair $(E, F)$ is naturally associated to a pair of paraballs. 

The reason for substituting a general quasiextremal pair $(E, F)$ by paraballs is that the parballs have nice geometric parametrizations. This is used in Section~\ref{S : mock distance} to define a mock-distance on the set of all paraballs similar to~\cite{Christ_extremal}. This in turn is used to prove that the images of distant paraballs under $X$ are essentially disjointly supported in $\Rd$.

In Section~\ref{S : proof of localization}, we combine the results from the previous two sections and complete the proof of Proposition~\ref{P : loc}.

\subsection{Notation}
$c, C$ will denote small and large positive implicit constants, respectively, whose precise value may change from one line to the next but is allowed to depend only on the dimension $d$ and the Lebesgue exponent $p$. The dependence of the implicit constants on other parameters will be denoted by additional subscripts. We write $A \lesssim B$ to denote $A \leq C B$, and write $A \sim B$ to denote $A \lesssim B \lesssim A$. For the remainder of the article we assume that $d \geq 3$ and all functions $f, g$ are nonnegative. For notational convenience, we denote $\|g\|_{\theta'} := \|g\|_{L_t^{q_\theta'} L_y^{r_\theta'}}$. 

\section{Proof of Theorem~\ref{T : main thm} : Existence and norm convergence}\label{S : Proof of main thm}
We state the following key propositions to be proved in subsequent sections. For $0 \leq \theta < 1$ with $(p_\theta, q_\theta, r_\theta)$ from~\eqref{E : def p_theta}, we write
\begin{equation}\label{E : def A_theta}
A_\theta := \|X\|_{L^{p_\theta} \to L_t^{q_\theta} L_y^{r_\theta}}, \qtq{and}  \scriptX(f, g) := \int_{\Rd} g \, Xf   
\end{equation}

\begin{proposition}{\bf(Frequency Localization)}\label{P : freq loc}
Let $\theta \in (0, 1)$. For each $\eps > 0$ there exist a positive $\delta$ and a finite $R$ such that the following holds. For each pair $(f, g)$ satisfying  $\scriptX(f, g) \geq (1 - \delta) A_\theta \|f\|_{L^{p_\theta}} \|g\|_{\theta'}$, there exist a symmetry $\phi^*$ so that with $A := \{|\phi^* f| < R \|f\|_{L^{p_\theta}}\}$ and $B := \{|\psi^* g| < R \|g\|_{\theta'}\}$, it holds that
$$
\scriptX(\rchi_A \phi^* f, \rchi_B \psi^* g) \geq (1 - \eps) A_\theta \|f\|_{L^{p_\theta}} \|g\|_{\theta'}.
$$
\end{proposition}

\begin{proposition}{\bf(Uniform localization)}\label{P : loc}
Let $\theta \in (0, 1)$. For each $\eps > 0$ there exist $\delta$ positive and a finite $R$ such that the following holds. For each $f$ satisfying $\|Xf\|_{L_t^{q_\theta} L_y^{r_\theta}} \geq A_\theta(1 - \delta) \|f\|_{L^{p_\theta}}$, there exists a symmetry $\phi^*$ so that
$$
\|\phi^* f\|_{L^{p_\theta} \big(\{|(s, x)| > R\} \cup \{|\phi^* f| > R \|f\|_{L^{p_\theta}}\}\big)} < \eps \|f\|_{L^{p_\theta}}.
$$
Moreover, the symmetry $\phi^*$ can be chosen to be independent of $\eps$.
\end{proposition}

Assuming the above propositions, we give the short proof of Theorem~\ref{T : main thm}.

\begin{proof}[Proof of Theorem~\ref{T : main thm}] The argument follows known patterns. By Fubini, at $\theta =  0$, each non-negative $f$ extremizes $X : L^1 \to L_t^{\infty} L_y^1$. Let $\{(f_n, g_n)\}$ be a normalized extremizing-pair of $X : L^{p_\theta} \to L_t^{q_\theta} L_y^{r_\theta}$, i.e., $\lim_n \scriptX(f_n, g_n) = A_\theta$ for some $\theta \in (0, 1)$. Being bounded, $f_n \rightharpoonup f, g_n \rightharpoonup g$ weakly in $L^{p_\theta}$ and $L_t^{q_\theta'}L_y^{r_\theta'}$ respectively. We assume that these, and all other limits arising in the proof, exist along a subsequence.

Since for a fixed smooth $\eta$ with compact support, the operator $G \mapsto \eta X^*(\eta G)$ maps $L^2$ to the Sobolev space $H^s$ for some $s > 0$ (~\cite[Theorem $8.11$]{CNSW}), for each fixed $R$ positive, the sequence $\scriptX(f_n^R, g_n^R) \to \scriptX(f^R, g^R)$ as $n \to \infty$ where $\Phi^R := \Phi \rchi_{\{|z| < R\} \cap \{|\Phi| < R\}}$ denotes the truncated and localized part of $\Phi$.

First, applying Proposition~\ref{P : loc} and then Proposition~\ref{P : freq loc} (and replacing $f_n$ by $\phi_n^* f_n$ and $g_n$ by $\psi_n^* g_n$), for $\tilde R > R$,
\begin{align*}
& A_\theta - o_R(1) \leq \lim_n \scriptX(f_n^R, g_n^{\tilde R}) + \scriptX(f_n^R, g_n \rchi_{\{|(t, y)| > \tilde R\}}) = \lim_n \scriptX(f_n^R, g_n^{\tilde R}) + o_{\tilde R}(1).
\end{align*}
Choosing $\tilde R$ large enough and replacing $R$ with $\tilde R$, we get $\scriptX(f^R, g^R) \geq A_\theta - o_R(1)$. Since $X$ is a positive operator, letting $R \to \infty$ implies $\scriptX(f, g) = A_\theta$. Furthermore,
$$
A_\theta = \scriptX(f, g) \leq A_\theta \|f\|_{L^{p_\theta}} \|g\|_{\theta'} \leq A_\theta.
$$
Thus, $\|f\|_{L^{p_\theta}} = 1$. By strict convexity~\cite[Theorem $2.11$]{LL}, we may now upgrade the weak convergence $f_n \rightharpoonup f$ to convergence in $L^{p_\theta}$, completing the proof.
\end{proof}

\section{Preliminaries}\label{S : prelims}
We start with a few preparatory results to be used repeatedly throughout this article. The following proposition roughly says : For $\theta \neq \theta_0$ in~\eqref{E : def theta_0} and $g$ with constant $t$-fiber (see condition~\eqref{E : def const $t$-fiber}), in order for the pair $(\rchi_E,g)$ to be large quasiextremal-pair of $X : L^{p_\theta} \to L_t^{q_\theta} L_y^{r_\theta}$, the function $g$ is necessarily ``almost" constant and $(\rchi_E, g)$ is a large quasiextremal-pair of $X : L^{p_{\theta_0}} \to L^{q_{\theta_0}}$.

For a precise statement, we need a bit of notation. We say that $g \in L_t^a L_y^b$ defined on $\R \times \Rk$ has constant $t$-fiber $B$ if there exists a measurable $T \subset \R$ so that
\begin{equation}\label{E : def const $t$-fiber}
B \rchi_{T}(t) \leq \|g(t, \cdot)\|_{L_y^b}^b \leq 2 B \rchi_{T}(t)    
\end{equation}
Note that for all such $g$ the $L_t^a L_y^b$ norm decouples : For each fixed $t_0 \in T$,
$$
\|g\|_{L_t^a L_y^b} \sim |T|^{1/a} \|g(t_0, \cdot)\|_{L_y^b}.
$$
The following is~\cite[Lemma $6.3$]{BS4}. As this will be used a few times, for the convenience of the reader, we note it below.
\begin{proposition}\label{P : intermediate strong}
Let $1 < s_j, u_j, v_j < \infty$ for $j = 0, 1$ with $v_0 < v_1$, and $L$ be a positive operator that satisfies the restricted weak type bounds
$$
\langle L(\rchi_E), \rchi_F \rangle \leq M_j |E|^{\frac 1{s_j}} \|\rchi_{F}\|_{L_t^{u_j'} L_y^{v_j'}}, \quad j = 0, 1.
$$
For each $\theta \in (0, 1)$ there exists a fixed positive constant $C_\theta = C_\theta(\theta, L)$ so that the following holds for $(1/s_\theta, 1/u_\theta, 1/v_\theta) := (1 - \theta) (1/s_0, 1/u_0, 1/v_0) + \theta (1/s_1, 1/u_1, 1/v_1)$. Let $0 < \eps < 1$. Let $g = \sum_{k \in \Z} 2^k g_k$, with constant $t$-fiber $B$ in $L_t^{u_\theta'} L_y^{v_\theta'}$, where $g_k \sim \rchi_{F_k}$, and measurable $E \subset \Rd$ satisfy $\langle L(\rchi_E), g \rangle \geq \eps M_0^{1 - \theta} M_1^{\theta} |E|^{\frac 1{s_\theta}} \|g\|_{L_t^{u_\theta'} L_y^{v_\theta'}}$. Then, with the distinguished frequency 
$
k_0 := (1 - \frac{v_\theta'}{v_0'})^{- 1} \Big[\theta Log(\frac{M_1}{M_0} + \frac{M_0}{M_1}) + (\frac{u_\theta'}{u_0'} - \frac{v_\theta'}{v_0'}) Log(1 + B) +(\frac 1{s_\theta} - \frac 1{s_0}) Log(1 + |E|) + (1 - \frac{u_\theta'}{u_0'}) Log(1 + \|g\|_{L_t^{u_\theta'} L_y^{v_\theta'}})\Big]
$
it holds that
\begin{equation}\label{E : Log delta^-1 terms}
\langle L(\rchi_E), g \rangle \leq (1 + \eps^{100}) \langle L(\rchi_E), \sum_{|k - k_0| \leq C_\theta Log(1 + \eps^{- C_\theta})} 2^k g_k \rangle.
\end{equation}
Additionally, for each integer $k$ in the above sum and $j = 0, 1$,
\begin{equation}\label{E : reduce to finite F_k}
\begin{gathered}
\frac{\eps^{C_\theta}}{C_\theta} \leq \frac{\langle L(\rchi_E), F_k \rangle}{|E|^{1/s_j} \|\rchi_{F_k}\|_{L_t^{u_j'} L_y^{v_{j'}}}}, \quad
\frac{\eps^{C_\theta}}{C_\theta} \leq \frac{\|\rchi_{F_k}\|_{L_t^{u_j'} L_y^{v_j'}}}{B^{\beta} |E|^{\Gamma_j}\|g\|_{L_t^{u'} L_y^{v'}}^{u' d_j}} \leq C_\theta \eps^{- C_\theta},
\end{gathered}
\end{equation}
where with $\Gamma := 1/{v_0'} - 1/{v_1'}$, the exponents $\beta := \big(\frac{1}{u_0'v_1'} - \frac{1}{v_0'u_1'}\big) u'/\Gamma, \Gamma_j := (1/{s_1} - 1/{s_0})/{v_j' \Gamma}$, and $d_j :=  1/{u_j'} + (1/{u_0} - 1/{u_1})/{v_j' \Gamma}$.
\end{proposition}

We say the pair $(f, g)$ is an $\eps$-quasiextremal pair at $\theta$ (for brevity we will mostly omit the $\theta$, if it is clear from context) if 
\begin{equation}\label{E : def eps-quasi pair}
\scriptX(f, g) \geq \eps \|f\|_{L^{p_\theta}} \|g\|_{\theta'}.
\end{equation} 
Abusing notation, if $f = \rchi_E$ and $g = \rchi_F$ are characteristic functions of measurable sets, we will simply say that $(E, F)$ is an $\eps$-quasiextremal pair if the above holds.

The next proposition is key to our analysis and will be applied repeatedly. It roughly says : If the images of two sets under $X$ are ``$\eps$-large" on the same part of $\Rd$, then they are of comparable volume up to a factor of $\sim \eps^{- C}$.

\begin{proposition}\label{P : similar vol E_1, E_2}
Let $\theta \in (0, 1)$ and $\eps$ be positive. Let $E_1, E_2$ be measurable subsets of $\Rd$ with the property that there exists a measurable $F \subset \Rd$ with constant $t$-fiber in $L_t^{q_\theta'} L_y^{r_\theta'}$ in case $\theta \neq \theta_0$, and these satisfy 
$$
X(\rchi_{E_i}) \geq \eps |E_i|^{1/{p_\theta}} \|\rchi_{F}\|_{\theta'} |F|^{- 1}, \quad i = 0, 1
$$ 
on a subset $\tilde F \subset F$ with $|\tilde F| \geq \eps |F|$. Then $\eps^C |E_1| \lesssim |E_2| \lesssim \eps^{- C} |E_1|$.
\end{proposition}
\begin{proof}
The proof uses a combination of the method of refinement of Christ and extrapolation. We start with the following lemma. For a pair of measurable sets $(E, F)$ we write $\scriptX(E, F) := \int_F X(\rchi_E)$.
\begin{lemma}\label{L : theta_0 strong}
Let $0 < \theta < 1, \theta \neq \theta_0$ and $\eps$ be positive. Let $(E, F)$ be an $\eps$-quasiextremal pair and $\rchi_F$ has constant $t$-fiber in $L_t^{q_\theta'} L_y^{r_\theta'}$. Then it holds that
$$
|E|^{1/p_{\theta_0}} |F|^{1/{q_{\theta_0}'}} \lesssim \eps^{- C} |E|^{1/p_\theta} \|\rchi_F\|_{\theta'}.
$$
\end{lemma}
\begin{proof}
For $\theta_0 < \theta < 1$ we apply bound~\eqref{E : reduce to finite F_k} for
$(s_0, u_0, v_0) = (p_{\theta_0}, q_{\theta_0}, r_{\theta_0})$ and $(s_1, u_1, v_1) =(p_1, q_1, r_1)$ to get
$$
|E|^{1/{p_\theta}} \|\rchi_F\|_{\theta'} \gtrsim \scriptX(E, F) \gtrsim \eps^C |E|^{1/p_{\theta_0}} |F_k|^{1/{q_{\theta_0}'}}.
$$
Similarly, for $0 < \theta < \theta_0$ applying Proposition~\ref{P : intermediate strong} for $(s_0, u_0, v_0) = (p_{\theta^2}, q_{\theta^2}, r_{\theta^2})$ and $(s_1, u_1, v_1) = (p_{\theta_0}, q_{\theta_0}, r_{\theta_0})$ gives the desired result.
\end{proof}

With this lemma in place, we now use the method of refinement. We give the proof for $d = 2 D$ even, for the odd case being identical. Assuming $|E_1| \leq |E_2|$, we will show that $|E_2| \lesssim \eps^{- C} |E_1|$. Denote $\alpha_i := \eps |E_i|^{1/{p_\theta}} \|\rchi_{F}\|_{\theta'} |F|^{- 1}, i = 1, 2$, and the average number of incidences in $E_1$ by $\beta := \frac{\scriptX(E_1, \tilde F)}{2 |E_1|} \gtrsim \eps^2 |E_1|^{- 1/{p_\theta'}} \|\rchi_{F}\|_{\theta'}$. We claim that there exist a point $(s_0, x_0) \in E_1$ and a measurable set $\Omega \subset \Rd$ with $|\Omega| \gtrsim \alpha_1^{D - 1} \alpha_2 \beta^D$ such that for each $(s_1, t_1, \ldots, s_D, t_D) \in \Omega$,
\begin{itemize}
\item $|t_i - t_j|\gtrsim \beta, \quad i \neq j$;
\item $|s_ j - s_{j - 1}| \gtrsim \alpha_1, \quad j \leq D - 1, \qtq{and} |s_D - s_{D - 1}| \gtrsim \alpha_2$;
\item $\Phi((s_1, t_1, \ldots, s_D, t_D)) := (s_D, x_0 - \sum_{j = 1}^D (s_{j - 1} -s_j) \gamma_0(t_j)) \in E_2$.
\end{itemize}
The proof of this follows from by now the well known standard argument known as the method of refinement of Christ~\cite{Christ_1998} (this in turn was inspired by and related to the work of Bourgain~\cite{Bourgain1986, Bourgain1991}, Schlag~\cite{Schlag1997, Schlag1998, Schlag1998kakeya}, and Wolff~\cite{Wolff1995, Wolff1997} on Kakeya-type problems, and is connected to the iterative argument of~\cite{Christ1985groups}), so we skip the details (see, e.g., \cite[Lemma $7.2$]{CB},~\cite[Lemma $5.6$]{BS4}, and \cite[bound (43)]{Tao-Wirght2001}).

By~\cite[Proposition $4.1$]{BS4}, for the jacobian $J\Phi \gtrsim \prod_{i = 1}^D |s_i - s_{i - 1} \, |\prod_{1 \leq j \neq i \leq D} |t_i - t_j|^2$, so by Bezout's theorem and a bit of arithmetic,
$$
|E_2| \geq |\Phi(\Omega)| \gtrsim \alpha_1^{2 (D - 1)} \alpha_2^2 \beta^{2 D^2 - D}.
$$
Substituting the values of $\alpha_i, \beta$, the above is equivalent to
$$
\frac{|E_2|}{|E_1|} \lesssim \Bigg(\frac{|E_1|^{1/p_{\theta_0}}|F|^{1/{q_{\theta_0}'}}}{\eps |E_1|^{1/{p_\theta}} \|\rchi_{F}\|_{\theta'}}\Bigg)^{\frac {(d^2 + d) p_\theta}{2 (2 - p_\theta)}}.
$$
For $\theta = \theta_0$, the claim is now obvious; whereas for $\theta \neq \theta_0$, an application of Lemma~\ref{L : theta_0 strong} finishes the proof. 
\end{proof}

\section{Lorentz Space estimates}\label{S : lorentz bound}
In proving Proposition~\ref{P : freq loc} our first step is to establish Lorentz space estimates for $X$. This roughly says : For a near extremizer pair $(f, g)$ of $X$, decomposing the functions into dyadic pieces $f = \sum 2^j f_j, g = \sum 2^k g_k$ where both $f_j, g_k \sim 1$, we can discard all small dyadic pieces of $f, g$. The following proposition establishes this for $f$ in the special case when $g$ possesses constant $t$-fiber.

\begin{proposition}\label{P : large dyadic of f}
Let $\theta_0 < \theta < 1$ and $\delta, \eta > 0$ be fixed. Let $g$ with constant $t$-fiber in $L_t^{q_\theta'} L_y^{r_\theta'}$, and $f = \sum_j 2^j f_j, f_j \sim \rchi_{E_j}$ be such that for each $j$, it holds that $2^{j p_\theta} |E_j| \sim \eta \|f\|_{L^{p_\theta}}^{p_\theta}$ and the number of incidences $\scriptX(E_j, g) \sim \delta |E_j|^{\frac1{p_\theta}} \|g\|_{\theta'}$. Then 
\begin{equation}\label{E : eps-improv for (f, g^l)}
\scriptX(f, g) \lesssim (\delta \eta)^c \|f\|_{L^{p_\theta}} \|g\|_{\theta'}.
\end{equation}
\end{proposition}

\begin{proof}
We may assume that $\|f\|_{L^{p_\theta}} = \|g\|_{\theta'} = 1$. Let $g = \sum_k 2^k g_k, g_k \sim \rchi_{F_k}$ for pairwise disjoint measurable $F_k \subset \Rd$. We consider the partition : For $m \in \Z$, let
\begin{equation}\label{E : def F_k^m}
F_k^m := F_k \cap \pi^{- 1}({\{t : 2^m \leq \int \rchi_{F_k} (t, y) dy < 2^{m + 1}\}}), \qtq{where} \pi (t, y) := t.
\end{equation}
By construction, each $F_k^m$ satisfies the constant $t$-fiber condition~\eqref{E : def const $t$-fiber} with $B = 2^m$. Fixing $\delta_1$ positive, for each $j$, we let $\Sigma_j^{\delta_1}$ denote the set of pairs $(k, m) \in \Z^2$ satisfying
$$
\delta_1 \leq \frac{\scriptX(E_j, F_k^m)}{|E_j|^{1/{p_\theta}} \|\rchi_{F_k^m}\|_{\theta'}} < 2 \delta_1.
$$
Denoting the subset $F_{k, j}^m := F_k^m \cap \{X (\rchi_{E_j}) \geq (1/2) \scriptX(E_j, F_k^m) |F_k^m|^{- 1}\}$, we have $\scriptX(E_j, F_k^m) \sim \scriptX(E_j, F_{k, j}^m)$. In a moment, we will prove the following.
\begin{align}
& \label{E : disjoint Fs} \sum_{j : (k, m) \in \Sigma_j^{\delta_1}} |F_{k, j}^m| \lesssim Log(1 + \delta_1^{- C}) |F_k^m| \qtq {for each fixed} k, m;
\\\label{E : delta_1 > delta}
& \Sigma_j^{\delta_1} \neq \emptyset \qtq{for some $j$ only if} \delta_1 \gtrsim \delta^{C}.
\end{align}

Assuming the above claims for now, we proceed. By our hypotheses on $2^{j p_\theta} |E_j|$,
\begin{equation}\label{E : neg-eta}
\scriptX(f, g) \sim \delta \sum 2^j |E_j|^{1/{p_\theta}} \|g\|_{\theta'} \lesssim \delta \eta^{- 1 + 1/{p_\theta}}.
\end{equation}
Additionally, since for $\delta_1$ and $j$ fixed, $\{F_{k, j}^m\subset F_k^m\}_{k, m}$ is a mutually disjoint collection, using the constant $t$-fiber property~\eqref{E : def const $t$-fiber} of $g$, the fact that $q_\theta < r_\theta$ for $\theta_0 < \theta < 1$, and claim~\eqref{E : disjoint Fs}, by H\"older's inequality,
\begin{align}\label{E : pos-eta}
\scriptX(f, g) & \notag \sim \sum_j 2^j \scriptX(E_j, \sum_{(k, m) \in \Sigma_j^{\delta_1}} 2^k \rchi_{F_{k, j}^m}) \lesssim \sum_j 2^j |E_j|^{1/{p_\theta}} \|\sum_{(k, m) \in \Sigma_j^{\delta_1}} 2^k \rchi_{F_{k, j}^m}\|_{\theta'}
\\\notag
& \leq \big(\sum_j (2^j |E_j|^{1/{p_\theta}})^{q_\theta}\big)^{1/{q_\theta}} \Big[\sum_j \int \Big(\sum_{(k, m) \in \Sigma_j^{\delta_1}} 2^{k r_\theta'} \int \rchi_{F_{k, j}^m}(t, y) \, dy \Big)^{\frac{q_\theta'}{r_\theta'}} \, dt \Big]^{1/q_\theta'}
\\\notag
& \lesssim \eta^{1/p_\theta - 1/q_\theta} Log(1 + \delta_1^{- C})^{1/{q_\theta'}} \big(B^{\frac{q_\theta'}{r_\theta'} - 1} \sum_k 2^{k r_\theta'} |F_k|\big)^{1/q_\theta'}
\\
& \lesssim \eta^{1/p_\theta - 1/q_\theta} Log(1 + \delta_1^{- C})^{1/q_\theta'}.
\end{align}
Since $p_\theta < q_\theta$, interpolating this with bound~\eqref{E : neg-eta}, and then summing over the dyadic values of $\delta_1 \in [\delta^C/C, 1]$ gives us~\eqref{E : eps-improv for (f, g^l)}. It remains to prove the above claims.

We first establish claim~\eqref{E : delta_1 > delta}. This follows from applying the following lemma to a fixed $E_j$ and then summing over the dyadic values of $0 < \eta_1 \leq 1$.

\begin{lemma}\label{L : delta_1 > delta}
Let $\theta_0 < \theta < 1$ and $0 < \delta_1 < 1$ be fixed. Let $g = \sum 2^k g_k, g_k \sim \rchi_{F_k}$, be a function with constant $t$-fiber in $L_t^{q_\theta'} L_y^{r_\theta'}$, and $F_k^m$ be as above. Suppose that for each $k, m$, the norm $\|2^k \rchi_{F_k^m}\|_{\theta'} \sim \eta_1 \|g\|_{\theta'}$ and let $E \subset \Rd$ a measurable set so that the number of incidences $\scriptX(E, F_k^m) \sim \delta_1 |E|^{1/{p_\theta}} \|\rchi_{F_k^m}\|_{\theta'}$. Then it holds that 
$$
\scriptX(E, g) \lesssim (\eta_1 \delta_1)^c |E|^{1/{p_\theta}} \|g\|_{\theta'}.
$$
\end{lemma}
\begin{proof}
By Proposition~\ref{P : intermediate strong} only $\sim Log(\delta^{- C})$ of the terms $2^k g_k$ contribute most of the incidences between $E$ and $g$ where $\delta := \frac{\scriptX(E, g)}{|E|^{1/{p_\theta}} \|g\|_{\theta'}}$. More precisely, by bound~\eqref{E : Log delta^-1 terms} corresponding to $(s_0, u_0, v_0) = (p_{\theta_0}, q_{\theta_0}, r_{\theta_0})$ and $(s_1, u_1, v_1) =(p_1, q_1, r_1)$, we have $\scriptX(E, g) \lesssim \scriptX(E, \sum_{|k - k_0| \lesssim Log(1 + \delta^{- C})} 2^k \rchi_{F_k})$ for some $k_0 = k_0 (E, \theta, g) \in \Z$. Thus, it suffices to prove that for each fixed $k$ in the sum
\begin{equation}
\scriptX(E, F_k) \lesssim (\eta_1 \delta_1)^c |E|^{1/{p_\theta}} \|\rchi_{F_k}\|_{\theta'}.
\end{equation}
This follows from a direct application of~\cite[the bound following $(7.5)$]{BS4}.
\end{proof}

Next, we turn to the proof of claim~\eqref{E : disjoint Fs}. We use an argument from~\cite{Christ_extremal} which has since appeared in related contexts in \cite{CB, BS4, NL}. First, we note that for $(k, m) \in \Sigma_j^{\delta_1}$,
$$
|E_j|^{1/p_{\theta_0}} |F_{k, j}^m|^{1/q_{\theta_0}'} \gtrsim \scriptX(E_j, F_{k, j}^m) \sim \delta_1 |E_j|^{1/{p_\theta}} \|\rchi_{F_k^m}\|_{\theta'}.
$$
Thus, by Lemma~\ref{L : theta_0 strong} and a bit of arithmetic, $|F_{k, j}^m|\gtrsim \delta_1^C |F_k^m|$. It suffices to prove that there exists a $C$ (possibly quite large), so that for each $k, m$ fixed,
\begin{equation}\label{E : sum_J}
\sum_{j \in J} |F_{k, j}^m| \lesssim |F_k^m|
\end{equation}
where $J$ is a $C Log(1 + \delta_1^{- C})$-separated set of integers.

We will proceed by the method of contradiction. By Cauchy-Schwartz,
\begin{align*}
|F_k^m|^{- 1} (\sum_j |F_{k, j}^m|)^2 \leq \sum_j |F_{k, j}^m| + \sum_{j_1 \neq j_2} |F_{k, j_1}^m \cap F_{k, j_2}^m|.
\end{align*}
Thus, invalidity of~\eqref{E : sum_J} implies that
$$
\delta_1^C |J|^2 |F_k^m| \lesssim |J|^2 \max_{j_1 \neq j_2} |F_{k, j_1}^m \cap F_{k, j_2}^m|.
$$
This guarantees existence of indices $j_1 \neq j_2 \in J$ so that $|F_{k, j_1}^m \cap F_{k, j_2}^m| \gtrsim \delta_1^C |F_k^m|$. Recalling that $|F_{k, j}^m| \gtrsim \delta_1^C |F_k^m|$, we apply Proposition~\ref{P : similar vol E_1, E_2} to the tuple $(E_{j_1}, E_{j_2}, F_{k, j_1}^m \cap F_{k, j_2}^m)$ to obtain $\delta_1^C |E_{j_1}| \lesssim |E_{j_2}| \lesssim \delta_1^{- C} |E_{j_2}|$. Since $2^{j_i p_\theta} |E_{j_i}| \sim \eta,$ for $ i = 1, 2$, this is equivalent to $2^{|j_1 -j_2|} \lesssim \delta_1^{- C}$ but this contradicts the fact that $J$ is $\sim Log(1 + \delta_1^{-C})$-separated. This completes the proof.
\end{proof}

We also have an analogue of Proposition~\ref{P : large dyadic of f}, where $g \in L_t^{q_\theta'} L_y^{r_\theta'}$ is now arbitrary, while $f$ is a characteristic function of a measurable set. To state this, we need a bit of notation. For $g \in L_t^{q_\theta'}L_y^{r_\theta'}$ and $l \in \Z$, (maintaining consistency with $F_k^m$ from~\eqref{E : def F_k^m}), we let $g^l$ denote the part of $g$ with constant $t$-fiber $2^l$, i.e.
\begin{equation}\label{E : def g^l}
g^l := g \rchi_{\pi^{- 1} (\{t : 2^l \leq \|g(t, \cdot)\|_{L_y^{r_\theta'}}^{r_\theta'} < 2^{l + 1}\})}.  
\end{equation}

\begin{proposition}\label{P : large dyadic of g}
Let $\theta_0 < \theta < 1$ and $\delta, \eta$ be positive. Let $E \subset \Rd$ be measurable and $g \in L_t^{q_\theta'}L_y^{r_\theta'}$ be such that whenever $\|g^l\|_{\theta'} > 0$, the norm $\|g^l\|_{\theta'} \sim \eta$ and the number of incidences $\scriptX(E, g^l) \sim \delta |E|^{1/{p_\theta}} \|g^l\|_{\theta'}$. Then it holds that
\begin{equation}\label{E : eps-improv for (E, g)}
\scriptX(E, g) \lesssim (\delta \eta)^c |E|^{1/{p_\theta}} \|g\|_{\theta'}.    
\end{equation}
In addition, there exists a collection of measurable subsets $E_l \subset E$ such that $\scriptX(E, g^l) \sim \scriptX(E_l, g^l)$ with $\sum_l |E_l| \lesssim Log(1 + \delta^{- C}) |E|$.
\end{proposition}

\begin{proof}
This is established in the course of the proof of \cite[Proposition $7.1$]{BS4} (more precisely~\cite[the bound following $(7.5)$]{BS4}). As the argument is rather long, we will not repeat it here.
\end{proof}

With $f = \sum_j 2^j f_j, f_j \sim \rchi_{E_j}$, for a mutually disjoint collection $\{E_j\}$, recalling that $\|f\|_{L^{p,s}} :\sim (\sum_j (2^j |E_j|^{\frac 1p})^s)^{\frac 1s}$, we have the following Lorentz space estimate.

\begin{theorem}\label{T : L(p, s) Lorentz}
Let $0 < \theta < 1$. For each $p_\theta < s < q_\theta$, the operator $X$ maps $L^{p_\theta, s}$ into $L_t^{q_\theta} L_y^{r_\theta}$.
\end{theorem}
\begin{proof}
First, we consider the case $\theta_0 < \theta < 1$. Let $\sum_j (2^j |E_j|^{1/{p_\theta}})^s = 1$ and $\|g\|_{\theta'} = 1$. Let $g \sim \sum_k 2^k g_k, g_k \sim \rchi_{F_k}$ for pairwise disjoint $F_k$. Fixing $F_k$, let us consider the partition $\{F_{k, l}^m \subset F_k\}_{l, m \in \Z}$,
\begin{equation}\label{E : def F_{k, l}^m}
F_{k, l}^m : = F_k \cap \pi^{- 1}(\{t : \|g(t, \cdot)\|_{L_y^{r_\theta'}}^{r_\theta'} \sim 2^l\}) \cap \pi^{- 1}(\{t : \int\rchi_{F_k} (t, y) dy \sim 2^m\}).    
\end{equation}

Fixing positive $\delta, \eta$, for the remainder of the proof, consider the collection of indices $(j, k, l, m)$ satisfying
$$
(2^j |E_j|^{1/p_\theta})^s \sim \eta, \qtq{and} \scriptX(E_j,F_{k, l}^m) \sim \delta |E_j|^{1/p_\theta} \|\rchi_{F_{k, l}^m}\|_{\theta'}.
$$
Fixing one $F_{k, l}^m$, we now run the exact argument leading to the proof of bound~\eqref{E : disjoint Fs}. This yields a collection of measurable subsets $F_{k, l, j}^m \subset F_{k, l}^m$ so that
$$
\scriptX(E_j, F_{k, l}^m) \sim \scriptX(E_j, F_{k, l, j}^m) \qtq{and} \sum |F_{k, l, j}^m| \lesssim Log(1 + \delta^{- C}) |F_{k, l}^m|.
$$
Thus, fixing $j$, replacing $F_{k, l}^m$ by $F_{k, l, j}^m$ recovers most incidences between $E_j$ and $g^l$. More precisely, with $g_j^l := \sum_k 2^k g_k \sum_m \rchi_{F_{k, l, j}^m}$, we have $\scriptX(E_j, g^l) \sim \scriptX(E_j, g_j^l)$, and additionally, arguing as in~\eqref{E : pos-eta} using that $q_\theta < r_\theta$ for $\theta_0 < \theta < 1$, we get
$$
\sum_j\|g_j^l\|_{\theta'}^{q_\theta'} \lesssim Log(1 + \delta^{- C}) \|g^l\|_{\theta'}^{q_\theta'}.
$$

Since the functions $g^l$ have disjoint support in $t$, we get $\sum_l \|g^l\|_{\theta'}^{q_\theta'} = \|g\|_{\theta'}^{q_\theta'} = 1$. Additionally, for a fixed $j$, for the same reason $\sum_l \|g_j^l\|_{\theta'}^{q_\theta'} = \|\sum_l g_j^l\|_{\theta'}^{q_\theta'}$. This lets us sum up to gain a positive power of $\eta$ : By H\"older's inequality,
\begin{align*}
\sum_j 2^j & \scriptX(E_j, g) \sim \sum_j 2^j \scriptX(E_j, \sum_l g_j^l) \lesssim \sum_j 2^j |E_j|^{1/{p_\theta}}\|\sum_l g_j^l\|_{\theta'}
\\
& \lesssim (\sum_j (2^j |E_j|^{1/{p_\theta}})^{q_\theta})^{1/{q_\theta}}(\sum_j \|\sum_l g_j^l\|_{\theta'}^{q_\theta'})^{1/{q_\theta'}} \lesssim \eta^{1/s - 1/{q_\theta}} Log(1 + \delta^{- C})^{1/{q_\theta'}}.
\end{align*}

On the other hand, by Lemma~\ref{L : delta_1 > delta}, each $\scriptX(E_j, g^l) \lesssim \delta^c |E_j|^{1/p_\theta} \|g^l\|_{\theta'}$ and this in turn by Proposition~\ref{P : large dyadic of g} implies that each $\scriptX(E_j, g) \lesssim \delta^c |E_j|^{\frac{1}{p_\theta}}\|g\|_{\theta'}$. In conclusion,
$$
\sum_j 2^j \scriptX(E_j, g) \lesssim \delta^c \eta^{\frac 1s - 1}.
$$
Interpolating the above two bounds and now summing over the dyadic values of $\eta, \delta \in (0, 1]$ completes the proof for $\theta_0 < \theta < 1$. For $0 < \theta \leq \theta_0$, we choose a $\theta_0 < \tilde \theta < 1$ and apply real interpolation to the trivial $L^1\rightarrow L_t^{\infty}L_y^1$ and the Lorentz space bound at $\tilde \theta$ (\cite[section $1.18.4-1.18.6$]{triebel}).
\end{proof}

For $g = \sum_lg^l$ from~\eqref{E : def g^l}, with $\|g\|_{L_t^{q,s} L_y^r} := (\sum_l \|g^l\|_{L_t^q L_y^r}^s)^{1/s}$ similar arguments give the following.

\begin{theorem}\label{T : L(q, s) lorentz}
Let $0 < \theta < 1$ and $p_\theta < s < q_\theta$. At $\theta = \theta_0$, the operator $X$ maps $L^{p_{\theta_0}}$ into $L^{q_{\theta_0}, s}$, and at $\theta \neq \theta_0$, it maps $L^{p_\theta}$ into $L_t^{q_\theta, s} L_y^{r_\theta}$.
\end{theorem}

\begin{proof}
At $\theta = \theta_0$, the proof of Theorem~\ref{T : L(p, s) Lorentz} works equally well for the adjoint $X^* : L^{q_{\theta_0}'} \to L^{p_{\theta_0}'}$, so we skip the repetitive details. Now, to prove it for $\theta \neq \theta_0$, by real interpolation, it suffices to consider when $\theta_0 < \theta < 1$. As above, it suffices to show the following : For $\delta, \eta$ fixed positive numbers, if $L^{p_\theta}$-normalized $f = \sum_j 2^j f_j, f_j \sim \rchi_{E_j}$ and $L_t^{q_\theta', s'} L_y^{r_\theta'}$-normalized $g = \sum g^l$ satisfy that
$$
\|g^l\|_{\theta'}^{s'}  \sim \eta, \qtq{and} \scriptX (E_j, g^l) \sim \delta |E_j|^{1/p_\theta} \|g^l\|_{\theta'},
$$
then the following $\eps$-improvement for the total number of incidences holds,
\begin{equation}\label{E : eps-improv for theta large}
\scriptX (f, g) \lesssim (\delta \eta)^c.
\end{equation}

We note that this is an upgrade from $f$ in the statement of Proposition~\ref{P : large dyadic of g} to general $f \in L^{p_\theta}$. Using the hypothesized bounds on $\|g^l\|_{\theta'}$, and summing over the dyadic values of $\eta \in (0, 1]$ in~\eqref{E : eps-improv for (f, g^l)}, we get $\scriptX(f, g^l) \lesssim \delta^c \|g^l\|_{\theta'}$. In consequence,
$$
\scriptX (f, g) \sim \sum_l \scriptX(f, g^l) \lesssim \delta^c \eta^{- 1/s}.
$$
The desired bound~\eqref{E : eps-improv for theta large} would follow from interpolating the above estimate with
$$
\scriptX (f, g) \lesssim \eta^c Log (1 + \delta^{- C}).    
$$

To obtain this estimate, fixing $E_j$, we apply Proposition~\ref{P : large dyadic of g} to get a collection of measurable subsets $E_{j, l} \subset E_j$ so that using H\"older's inequality the following holds.
\begin{align*}
\scriptX (f, g) & \sim \sum_l \scriptX(\sum_j 2^j \rchi_{E_{j, l}}, g^l) \lesssim \sum_l \|g^l\|_{\theta'} \|\sum_j 2^j \rchi_{E_{j, l}}\|_{L^{p_\theta}}
\\
& \leq (\sum_l \|\sum_j 2^j \rchi_{E_{j, l}}\|_{L^{p_\theta}}^{p_\theta})^{1/p_\theta} (\sum_l \|g^l\|_{\theta'}^{p_\theta'})^{1/p_\theta'} \lesssim \eta^{\frac 1{s'} (1 - \frac{q_\theta'}{p_\theta'})} Log (1 + \delta^{- C})^{1/p_\theta}.
\end{align*}
Since $p_\theta < q_\theta$, this finishes the proof. 
\end{proof}

\section{Proof of Proposition~\ref{P : freq loc} : Frequency Localization}\label{S : proof of fre loc}

We now have the necessary tools required to prove Proposition~\ref{P : freq loc}.

\begin{proof}[Proof of Proposition~\ref{P : freq loc}] 
We first consider the simpler case, $\theta = \theta_0$, for which we are looking at the operator $X : L^{p_{\theta_0}} \to L^{q_{\theta_0}}$. The proof for this case follows an argument of Christ in~\cite[Lemma $6.1$]{Christ_extremal}. Let $f = \sum_j 2^j f_j, f_j \sim \rchi_{E_j}$, and $g = \sum_k 2^k g_k, g_k \sim \rchi_{F_k}$, for pairwise disjoint collections of measurable sets $\{E_j\}, \{F_k\}$. We may assume that both $f, g$ are normalized. Thus, $\sum 2^{j p_{\theta_0}} |E_j| \sim 1 \sim \sum 2^{k q_{\theta_0}'} |F_k|$. Suppose that $\scriptX(f, g) \geq(1 - \delta) A_{\theta_0}$ for a sufficiently small positive $\delta$.

We start by truncating $f, g$ by deleting their $\delta$-small part. More precisely, we consider $S$ the set of all $j \in \Z$ so that $2^{j p_{\theta_0}} |E_j| > \delta$. We construct a similar index set $T$ associated with $g$. Finally, we consider the truncated functions
$$
\tilde f := \sum_{j \in S} 2^j f_j \qtq{and} \tilde g := \sum_{k \in T} 2^k g_k.
$$
By Theorems~\ref{T : L(p, s) Lorentz} and~\ref{T : L(q, s) lorentz} and the triangle inequality, $\scriptX(\tilde f, \tilde g) \geq(1 - C \delta^c) A_{\theta_0}$.

\begin{lemma}\label{L : freq bound theta_0}
There exist integers $J_0, K_0$ so that
$$
|j - J_0| \lesssim \delta^{- C}, \quad  j \in S \qtq{and} |k - K_0| \lesssim \delta^{- C}, \quad k \in T.
$$
\end{lemma}

At $\theta = \theta_0$, the conclusion of the proposition is now immediate modulo the proof of the above lemma. Indeed, given any $\eps > 0$ we may now take $\delta := c \eps^C$ and $R := C \eps^{- C}$ for some fixed constants $c, C$, and choose an appropriate scaling symmetry $\phi = S_{(\alpha, \beta)}$ where 
$$
2^{- J_0 p_{\theta_0}} = \alpha^d \beta^{\frac{d^2 - d}2} \qtq{and}  2^{- K_0 q_{\theta_0}'} = \alpha^{d - 1} \beta^{\frac{d^2 - d}2 + 1}.
$$

\begin{proof}[Proof of Lemma~\ref{L : freq bound theta_0}]
Throughout what follows, we fix $0 < \delta < 1$. We first show that $S$ lies inside an interval of length~$\sim \delta^{- C}$. With $D := \max_{j_1, j_2 \in S} |j_1 - j_2|$ and $N := |S| + |T|$, we decompose $S = S_L \sqcup S_R$ so that $|j_1 -j_2| \geq D/N$ for each $(j_1, j_2) \in S_L \times S_R$. With this partition in place, fixing a large $\lambda > 0$ to be chosen precisely momentarily, we partition each
$F_k$ into $F_k =:  F_{k L} \sqcup F_{k R} \sqcup F_{k s}$
where
\begin{align*}
F_{k L} :=\{(t, y) \in F_k : X \rchi_{E_j} (t, y) \geq \delta^\lambda |E_j|^{1/p_{\theta_0}} |F_k|^{- 1/q_{\theta_0}} \, \text{for some $j \in S_L$}\},
\\
\text{and} \, F_{k R} :=\{(t, y) \in F_k : X \rchi_{E_j} (t, y) \geq \delta^\lambda |E_j|^{1/p_{\theta_0}} |F_k|^{- 1/q_{\theta_0}} \, \text{for some $j \in S_R$}\} \setminus F_{k L}.
\end{align*}
Finally, we perform a decomposition of $\tilde f$ into $\tilde f = f_L + f_R$ where
$$
f_L := \sum_{j \in S_L} 2^j f_j, \qtq{and} f_R := \sum_{j \in S_R} 2^j f_j,
$$ 
and decompose $\tilde g$ into disjointly supported pieces, $\tilde g = g_L + g_R + g_s$ where
\begin{equation}\label{E : partition g}
g_L := \sum 2^k g_k \rchi_{F_{k L}}, \, g_R := \sum 2^k g_k \rchi_{F_{k R}},\qtq{and} g_s : = \sum 2^k g_k \rchi_{F_{k s}}.
\end{equation}

By construction, for each $j \in S$ and $k \in T$, the number of incidences
$$
\scriptX(2^j f_j, 2^k \rchi_{F_{k s}}) \lesssim \delta^{\lambda} 2^{j + k} |E_j|^{1/p_{\theta_0}} |F_k|^{1/q_{\theta_0}'} \lesssim \delta^{\lambda}.
$$
Summing over $N^2 \lesssim \delta^{- C}$ many such pairs, we get $\scriptX(\tilde f, g_s) \lesssim \delta^{\lambda - C}$. Additionally, applying H\"olders inequality with $p = p_{\theta_0}, q = q_{\theta_0}, \eps_0 :=  1 - p/q$,
\begin{align}\label{E : small pure incidences}
\sum_{i = L, R} \scriptX(f_i, g_i) \leq  A_{\theta_0} \sum_{i = L, R} &\notag \|f_i\|_{L^p} \|g_i\|_{L^{q'}} \leq  A_{\theta_0} (\sum_{i = L, R} \|f_i\|_{L^p}^q)^{1/q} (\sum_{i = L, R} \|g_i\|_{L^{q'}}^{q'})^{1/{q'}}
\\
& \leq A_{\theta_0} \max(\|f_L\|_{L^p}, \|f_R\|_{L^p})^{\eps_0} \leq A_{\theta_0} (1 - \delta)^{\eps_0}.
\end{align}
For $\lambda$ large enough this implies that the number of cross incidences
\begin{equation}\label{E : large cross insidences}
\scriptX(f_L, g_R) + \scriptX(f_R, g_L) \gtrsim \delta^C.
\end{equation}

We address when $\scriptX(f_L, g_R) \gtrsim \delta^C$, the other case being identical. Since $N^2 \lesssim \delta^{- C}$, by the triangle inequality we get a $\sim \delta^C$-quasiextremal pair $(E_j, F_{k R})$ with $j \in S_L$. Indeed, fixing a positive large parameter $\lambda$ to be chosen momentarily, suppose that none of the pairs $(E_j, F_{k R})$ is a $\delta^\lambda$-quasiextremal pair. As a result,
\begin{align}\label{E : get large quasiextremal}
\scriptX(f_L, g_R) \lesssim \scriptX (\sum 2^j f_j, \sum 2^k \rchi_{F_{k R}}) \lesssim \delta^{\lambda} \sum 2^j |E_j|^{1/{p_{\theta_0}}} 2^k |F_{k R}|^{1/{q_{\theta_0}'}} \lesssim \delta^{\lambda  - C},
\end{align}
producing the desired conclusion for $\lambda$ sufficiently large. Once more, by the triangle inequality and some arithmetic, we get that 
$$
X \rchi_{E_j} \geq \delta^C |E_j|^{1/p_{\theta_0}} |F_{k R}|^{- 1/q_{\theta_0}}
$$
on a measurable subset $F_1 \subset F_{k R}$ satisfying $|F| \gtrsim |F_{k R}|$.

On the other hand, since $|S| \lesssim \delta^{- C}$, we may assume that for some $j' \in S_R$
$$
X \rchi_{E_{j'}} \geq \delta^C |E_{j'}|^{1/p_{\theta_0}} |F_{k R}|^{- 1/q_{\theta_0}}
$$
on a measurable further subset $F \subset F_1$ satisfying $|F| \gtrsim \delta^C |F_{k R}|$. Applying Proposition~\ref{P : similar vol E_1, E_2} we deduce that $\delta^C |E_j| \lesssim |E_{j'}| \lesssim \delta^{-C} |E_j|$. Combined with the assumption that each $\delta \leq 2^{l p_{\theta_0}} |E_l| \leq 1$, we get $|j-j'| \lesssim Log(\delta^{- C})$. Since the pair $(j, j') \in S_L \times S_R$, this implies that $|j_1 - j_2| \lesssim \delta^{- C}$ for all $j_1, j_2 \in S$.

The above argument works equally well for the adjoint $X^* : L^{q_{\theta_0}'} \to L^{p_{\theta_0}'}$, giving the desired estimate on $T$. This completes the proof of the lemma.
\end{proof}

We now consider the remaining cases, $\theta \neq \theta_0$, corresponding to the mixed norm bounds $X : L^{p_\theta} \to L_t^{q_\theta} L_y^{r_\theta}$. For this we extend the method in~\cite[Lemma $6.1$]{Christ_extremal} used in the above case to a mixed norm setting. One of the obstructions stems from the fact that the norm of a function in a mixed norm space does not necessarily decouple if it is decomposed into disjointly supported pieces, which was used to arrive at the bound~\eqref{E : small pure incidences}. Thus, for example, a direct analogue of~\eqref{E : partition g} might not result in a nontrivial number of cross-incidences as in bound~\eqref{E : large cross insidences}. To overcome this, we proceed as follows.

Let $f = \sum_j 2^j f_j, f_j \sim \rchi_{E_j}$, and $g = \sum_k 2^k g_k, g_k \sim \rchi_{F_k}$, be normalized. Thus, $\sum 2^{j p_\theta} |E_j| \sim 1$, and recalling $g^l$ from~\eqref{E : def g^l}, we have $\sum \|g^l\|_{\theta'}^{q_\theta'} = 1$. Suppose that $\scriptX(f, g) \geq(1 - \delta) A_\theta$ for a sufficiently small positive $\delta$.

First, we truncate $f$. Let $S := \{j \in \Z : 2^{j p_{\theta}} |E_j| > \delta\}$. By Theorem~\ref{T : L(p, s) Lorentz} and the triangle inequality, for $F := \sum_{j \in S} 2^j f_j$, we have $\scriptX(F, g) \geq (1 - C \delta^c) A_\theta$.

Next, we truncate $g$. Let $L$ be the set of all $l \in \Z$ so that $\|g^l\|_{\theta'}^{q_\theta'} > \delta$ and $(\rchi_{E_j}, g^l)$ is a $\delta^C$-quasiextremal pair for some $j = j(l) \in S$. Once again, for a sufficiently large $C$, by Theorem~\ref{T : L(q, s) lorentz} and arguing as in~\eqref{E : get large quasiextremal}, we get $\scriptX(F, \sum_L g^l) \geq (1 - C \delta^c) A_\theta$. Now, fixing $l \in L$, using the hypothesized bound on $\scriptX(E_{j (l)}, g^l)$, we apply Proposition~\ref{P : intermediate strong}. This gives us a distinguished frequency $K(j, l) \in \Z$ so that
\begin{equation}\label{E : few l}
\scriptX(F, \sum_{l \in L} \sum_{k \in \scriptK_l} 2^k g_k \rchi_{F_{k, l}}) \geq (1 - C \delta^c) A_\theta,
\end{equation}
where $F_{k, l} : = F_k \cap \pi^{- 1}(\{t : 2^l \leq \|g(t, \cdot)\|_{L_y^{r_\theta'}}^{r_\theta'} < 2^{l + 1}\})$ and 
$$
\scriptK_l := \bigcup \{k \in \Z : |k - K(j, l)| < C Log (1 + \delta^{- C})\}.
$$
The above union is taken over all $j \in S$ such that $(\rchi_{E_j}, g^l)$ is a $\delta^C$-quasiextremal pair. For the final step of the truncation of $g$, fixing $l \in L$ and $k \in \scriptK_l$, let $\scriptM(l, k)$ denote the set of all $m \in \Z$ such that $F_{k, l}^m$ in~\eqref{E : def F_{k, l}^m} satisfies that $2^k \|\rchi_{F_{k, l}^m}\|_{\theta'} > \delta^C$. Once again, for a sufficiently large value of $C$, by Theorem~\ref{T : L(q, s) lorentz} and the triangle inequality, for the truncated function
$$
G := \sum_{l \in L} \sum_{k \in \scriptK_l} \sum_{m \in \scriptM(l, k)} 2^k g_k \rchi_{F_{k, l}^m},
$$
from bound~\eqref{E : few l} we have
\begin{equation}\label{E : large X(F, G)}
\scriptX(F, G) \geq(1 - C \delta^c) A_{\theta}.
\end{equation}
Furthermore, we recall that for the distinguished frequencies $K(j, l)$ we have 
\begin{equation}\label{E : nearby k(j, l)}
|K(j, l) - K(j', l')| = |C_1 (j - j') + C_2 (l - l')| + \scriptO(Log(1 + \delta^{- C})).
\end{equation}
for some fixed constants $C_1, C_2$. We now prove that the following holds.

\begin{lemma}\label{L : freq bound non theta_0}
There exist integers $j_0, l_0$ such that for each $j \in S$ and $l \in L$
$$
|j - j_0| + |l - l_0| \lesssim \delta^{- C}.
$$
\end{lemma}

As in the case of $\theta = \theta_0$, combining bounds~\eqref{E : large X(F, G)} and~\eqref{E : nearby k(j, l)}, we may now conclude the proof of Proposition~\ref{P : freq loc} modulo the proof of Lemma~\ref{L : freq bound non theta_0} by choosing an appropriate scaling symmetry $\phi = S_{(\alpha, \beta)}$ where 
$$
2^{- j_0 p_\theta} = \alpha^d \beta^{\frac{d^2 - d}2} \qtq{and}  2^{- K(j_0, l_0)} = \beta^{- \frac 1{q_\theta'}} (\alpha^{d - 1} \beta^{\frac{d^2 - d}2})^{- \frac 1{r_\theta'}}.
$$
\end{proof}

\begin{proof}[Proof of Lemma~\ref{L : freq bound non theta_0}]
For the remainder of the proof, we fix $0 < \delta < 1$. We first establish that $S$ is contained in an interval of length $\sim \delta^{- C}$. For $R > 0$ and a pair of measurable sets $(E, F)$, we say a point $(t, y) \in F$ is $R$-strong if
\begin{equation}\label{E : def lambda-large}
X \rchi_E (t, y) \geq \delta^R |E|^{1/p_\theta} \|\rchi_F\|_{\theta'} |F|^{- 1}.
\end{equation}
We let $R(E, F)$ denote the set of all $R$-strong points of $(E, F)$. We note that by Theorem~\ref{T : L(q, s) lorentz} and a bit of arithmetic, for any positive $R$ and given a $2 \delta^R$-quasiextremal pair $(E, F)$, the set of all $R$-strong points of the pair $(E, F)$ is a large subset of $F$; more precisely, $|R(E, F)| \gtrsim \delta^{C R} |F|$.

With $M := \max_{j_1, j_2 \in S} |j_1 - j_2|$ and $N := |S| + |L| \lesssim \delta^{- C}$, we consider a decomposition $S = S_U \sqcup S_D$ so that $|j_1 -j_2| \geq M/N$ for each $(j_1, j_2) \in S_U \times S_D$. With this partition in place, fixing a large $\lambda > 0$ to be chosen precisely momentarily, we recall the decomposition of $F_k$ into $F_{k, l}^m$ in~\eqref{E : def F_{k, l}^m}, and now further decompose each 
$$
F_{k, l}^m =: F_{k, l}^{m, U} \sqcup F_{k, l}^{m, D} \sqcup F_{k, l}^{m, s}
$$
where (we recall that $\rchi_{F_{k, l}^m}$ satisfies the constant $t$-fiber property~\eqref{E : def const $t$-fiber} with $B = 2^m$)
\begin{align}\label{E : def large fiber}
& F_{k, l}^{m, U} := \bigcup_{j \in S_U} \lambda (E_j, F_{k, l}^m) \cap \pi^{- 1}\big(\{t : |\{y : (t, y) \in  \lambda (E_j, F_{k, l}^m)\}| \geq \delta^\lambda 2^m\}\big)
\\\notag
& F_{k, l}^{m, D} := \bigcup_{j \in S_D} \lambda (E_j, F_{k, l}^m) \cap \pi^{- 1}\big(\{t : |\{y : (t, y) \in  \lambda (E_j, F_{k, l}^m)\}| \geq \delta^\lambda 2^m\}\big) \setminus F_{k, l}^{m, U}.
\end{align}
Finally, we perform a decomposition of $F$ into $F = f_U + f_D$ where
$$
f_U := \sum_{j \in S_U} 2^j f_j, \qtq{and} f_D := \sum_{j \in S_D} 2^j f_j
$$ 
and consider a corresponding decomposition of $G = g_U + g_D + g_s$ where
$$
g_U := \sum 2^k g_k \rchi_{F_{k, l}^{m, U}}, \, g_D := \sum 2^k g_k \rchi_{F_{k, l}^{m, D}},\qtq{and} g_s : = \sum 2^k g_k \rchi_{F_{k, l}^{m, s}}.
$$

By construction, for each $j \in S, l \in L, k \in \scriptK_l, m \in \scriptM(l, k)$, the number of incidences
\begin{equation}\label{E : small incidences}
\scriptX(2^j f_j, 2^k g_k \rchi_{F_{k, l}^{m, s}}) \lesssim \delta^{\lambda} 2^{j + k} |E_j|^{\frac{1}{p_\theta}} \|\rchi_{F_{k, l}^m}\|_{\theta'} \lesssim \delta^{\lambda}.
\end{equation}
As before, summing over $N^4 \lesssim \delta^{- C}$ such tuples, we get $\scriptX(F, g_s) \lesssim \delta^{\lambda - C}$. Next, we claim that choosing $\lambda$ large enough depending on $\theta$, yields
\begin{equation}\label{E : large incidence for f_D, g_U}
\scriptX(f_U, g_D) + \scriptX(f_D, g_U) \gtrsim \delta^C.
\end{equation}

To prove the above claim, we first consider the case $\theta_0 < \theta < 1$. In this region, $q_\theta' > r_\theta'$. Thus, by H\"older's inequality and mutually disjoint supports of $g_U, g_D$,
\begin{align*}
\scriptX(f_U, g_U) & + \scriptX(f_D, g_D) \leq A_\theta (\|f_U\|_{L^{p_\theta}}^{q_\theta} + \|f_D\|_{L^{p_\theta}}^{q_\theta})^{1/{q_\theta}} (\|g_U\|_{\theta'}^{q_\theta'} + \|g_D\|_{\theta'}^{q_\theta'})^{1/{q_\theta'}}
\\
& \leq A_\theta (1 - \delta)^{1 - \frac{p_\theta}{q_\theta}} (\int \max_{i = U, D} \|g_i (t, \cdot)\|_{L_y^{r_\theta'}}^{q_\theta' - r_\theta'} \|g (t, \cdot)\|_{L_y^{r_\theta'}}^{r_\theta'} dt)^{1/q_\theta'} \leq A_\theta (1 - \delta)^{1 - \frac{p_\theta}{q_\theta}}.
\end{align*}
Combining this with the fact that $\scriptX(F, g_s) \lesssim \delta^{\lambda - C}$, choosing $\lambda$ sufficiently large depending on $\theta$, we deduce claim~\eqref{E : large incidence for f_D, g_U}.

Next, we prove claim~\eqref{E : large incidence for f_D, g_U} for the remaining values of $\theta$ where $0 < \theta < \theta_0$. In this region, $q_\theta > r_\theta$. We proceed by the method of contradiction. Suppose, if possible, that $\scriptX(f_U, g_D) + \scriptX(f_D, g_U) < \delta^\lambda$. Thus, with $H_i := X f_i \rchi_{\supp(g_i)}$, by the triangle inequality in $L_t^{q_\theta/r_\theta}$ and the fact that $p_\theta < r_\theta$,
\begin{align*}
A_\theta^{r_\theta} (1 - C \delta^\lambda) \leq \Big[\int (\sum_{i = U, D} \int H_i(t, y)^{r_\theta} & dy)^{\frac{q_\theta}{r_\theta}} dt\Big]^{\frac{r_\theta}{q_\theta}} \leq \sum_{i = U, D}  \|H_i\|_{L_t^{q_\theta} L_y^{r_\theta}}^{r_\theta}
\\
& \leq A_\theta^{r_\theta} \sum_{i = U, D}  \|f_i\|_{L^{p_\theta}}^{r_\theta} \leq A_\theta^{r_\theta} \max_{i = U, D}  \|f_i\|_{L^{p_\theta}}^{r_\theta - p_\theta}.
\end{align*}
For $\lambda$ sufficiently large, this contradicts the fact that each $\|f_i\|_{L^{p_\theta}}^{p_\theta} \geq \delta$, completing the proof of claim~\eqref{E : large incidence for f_D, g_U}.

Looking at bound~\eqref{E : large incidence for f_D, g_U}, we address when $\scriptX(f_U, g_D) \gtrsim \delta^C$, the other case being identical. Since $N \lesssim \delta^{- C}$, we get a $\sim \delta^C$-quasiextremal pair $(E_j, F_{k, l}^{m, D})$ with $j \in S_U$ and additionally, arguing as in bound~\eqref{E : small incidences} using the fact that $\rchi_{F_{k, l}^m}$ has constant $t$-fiber, we get $|\pi(F_{k, l}^{m, D})| \gtrsim \delta^C |\pi(F_{k, l}^m)|$. Thus, again by the constant $t$-fiber property and by construction, $|F_{k, l}^{m, D}| \gtrsim \delta^C |F_{k, l}^m|$. As in the remark following~\eqref{E : def lambda-large}, the associated $C$-strong set $C(E_j, F_{k, l}^m)$ is of measure at least $\gtrsim \delta^C |F_{k, l}^m|$.

Furthermore, since $|S| \lesssim \delta^{- C}$ and $|F_{k, l}^{m, D}| \gtrsim \delta^C |F_{k, l}^m|$, for the intersection of $C$-strong points $|C (E_j, F_{k, l}^m) \cap C (E_{j'}, F_{k, l}^m)| \gtrsim \delta^C |F_{k, l}^m|$ for some $j' \in S_D$. We now argue as in the proof of Lemma~\ref{L : freq bound theta_0}, first apply Proposition~\ref{P : similar vol E_1, E_2} and then use the fact that each $2^{j p_\theta} |E_j| \geq \delta$, to deduce that $S$ lies in an interval of length $\sim \delta^{- C}$.

Next, we prove that $L$ lies within an interval of length $\sim \delta^{- C}$. We write the set $L = L_\alpha \sqcup L_\beta$ so that for $(l_1, l_2) \in L_\alpha \times L_\beta$ we have $|l_1 - l_2| \geq \frac{\max_{l, l'} |l - l'|}N$ where $N := |S| + |L| \lesssim \delta^{- C}$. Fixing a large positive parameter $\lambda$, to be chosen precisely momentarily, we first partition each $E_j =: E_j^\alpha \sqcup E_j^\beta \sqcup E_j^s$ where
\begin{align*}
& E_j^\alpha := \bigcup_{l \in L_\alpha} \{X^*G^l \geq \delta^\lambda |E_j|^{- 1/{p_\theta'}} \|G^l\|_{\theta'}\} \cap E_j,
\\
& E_j^\beta := \bigcup_{l \in L_\beta} \{X^*G^l \geq \delta^\lambda |E_j|^{- 1/{p_\theta'}} \|G^l\|_{\theta'}\} \cap E_j \setminus E_j^\alpha, 
\end{align*}
and now consider the associated decomposition $G = g_\alpha + g_\beta$ where
$$
g_\alpha : = \sum_{l \in L_\alpha} G^l, \qquad g_\beta = \sum_{l \in L_\beta} G^l
$$
and similarly write $F = f_\alpha + f_\beta + f_s$ where
$$
f_\alpha := \sum_j 2^j f_j \rchi_{E_j^\alpha}, \qquad f_\beta := \sum_j 2^j f_j \rchi_{E_j^\beta}, \qtq{and} f_s := \sum_j 2^j f_j \rchi_{E_j^s}
$$

By construction, $g_\alpha$ and $g_\beta$ have disjoint $t$-support, so $\|G\|_{\theta'}^{q_\theta'} = \|g_\alpha\|_{\theta'}^{q_\theta'} + \|g_\beta\|_{\theta'}^{q_\theta'}$. This allows us to repeat the argument leading to the bound~\eqref{E : large cross insidences}, and to conclude that choosing $\lambda$ large enough depending on $\theta$, we get $\scriptX(f_\alpha, g_\beta) + \scriptX(f_\beta, g_\alpha) \gtrsim \delta^C$. We consider when $\scriptX(f_\alpha, g_\beta) \gtrsim \delta^C$, identical arguments work for the other case.

As in~\eqref{E : get large quasiextremal}, we get a $\delta^C$-quasiextremal pair $(\rchi_{E_j^\alpha}, G^l)$ where $l \in L_\beta$. Next, by Proposition~\ref{P : intermediate strong} this in turn produces a pair $(E_j^\alpha, F_{k l})$ which is a $\sim \delta^C$-quasiextremal pair of $X : L^{p_{\theta_0}} \to L^{q_{\theta_0}}$ and 
\begin{equation}\label{E : large vol of F_kl}
\delta^C |E_j^\alpha|^{\frac{2 (d - 1) q_{\theta_0}'}{d^2 + d}} \lesssim 2^{q_{\theta _0}' r_\theta l} |F_{k l}| \lesssim \delta^{- C} |E_j^\alpha|^{\frac{2 (d - 1) q_{\theta_0}'}{d^2 + d}}.
\end{equation}
Thus, on a measurable subset $E \subset E_j^\alpha$ satisfying $|E| \gtrsim |E_j^\alpha|$, it holds that
\begin{equation}\label{E : large X*F1}
X^*(\rchi_{F_{k l}}) \gtrsim \delta^C |E_j^\alpha|^{- 1/{p_{\theta_0}'}} |F_{k l}|^{1/{q_{\theta_0}'}} \geq \delta^C |E|^{- 1/{p_{\theta_0}'}} |F_{k l}|^{1/{q_{\theta_0}'}}.
\end{equation}

On the other hand, since $L$ contains at most $\lesssim \delta^{- C}$ many points, for some $l' \in L_\alpha$
$$
X^*G^{l'}\geq\delta^C|E_j|^{- 1/{p_\theta'}}\|G^{l'}\|_{\theta'} \geq \delta^C|E|^{- 1/{p_\theta'}}\|G^{l'}\|_{\theta'}
$$
on a further subset $\tilde E \subset E$ satisfying $|\tilde E| \gtrsim \delta^C |E|$. Consequently, $(\rchi_{E} ,G^{l'})$ is a $\sim \delta^C$-quasiextremal pair. As above, applying Proposition~\ref{P : intermediate strong} once more we get $(E, F_{k' l'})$ a $\sim \delta^C$-quasiextremal pair of $X : L^{p_{\theta_0}} \to L^{q_{\theta_0}}$ with 
\begin{equation}\label{E : large vol of F_kl'}
\delta^C |E|^{\frac{2 (d - 1) q_{\theta_0}'}{d^2 + d}} \lesssim 2^{q_{\theta_0}' r_\theta l'} |F_{k' l'}| \lesssim \delta^{- C} |E|^{\frac{2 (d - 1) q_{\theta_0}'}{d^2 + d}}.
\end{equation}
Thus, on a subset $E_1 \subset E$ with $|E_1| \gtrsim |E|$, it holds that  
\begin{equation}\label{largeX*F2}
X^*(\rchi_{F_{k' l'}}) \gtrsim \delta^C |E|^{- 1/{p_{\theta_0}'}} |F_{k' l'}|^{1/{q_{\theta_0}'}}.
\end{equation}
Using the lower bounds~\eqref{E : large X*F1} and\eqref{largeX*F2}, the method of refinement as in the proof of Proposition~\ref{P : similar vol E_1, E_2} (see, e.g., \cite[bound $7.18$]{BS4}) implies that $\delta^C |F_{k l}| \lesssim |F_{k' l'}| \lesssim \delta^{- C} |F_{k l}|$. By the bounds~\eqref{E : large vol of F_kl} and~\eqref{E : large vol of F_kl'}, this in turn implies that $|l - l'| \lesssim Log(\delta^{- C})$. Thus, $|l_1 - l_2| \lesssim \delta^{- C}$ for all $l_1, l_2 \in  L$. Finally, the proof of Lemma~\ref{L : freq bound non theta_0} is complete.
\end{proof}

\section{Paraballs and Quasiextremal pairs}\label{S : paraballs}
So far, given $\eps$-quasiextremal pair $(f, g)$, we are able to extract a related $\eps^C$-quasiextremal pair $(E, F)$. However, there is a possibility that $E$ might turn out to be a disjoint union of arbitrarily large number of small subsets with arbitrarily large mutual distances (see, e.g., ~\cite[Remark $13.1$]{Christ_quasiextremal}). This is certainly not desirable for us. To show that this does not occur for the operator $X$, our next task is to identify a class of ``natural" quasiextremal pairs of $X$ and then relate the above pair $(E, F)$ to one belonging to this natural class (see Theorem~\ref{T : reduce to paraball}).

\begin{definition}\label{E : def paraball}
Set $B(0, 0, 0, 1, 1) = B^*(0, 0, 0, 1, 1)$ to be the unit cube $\{z \in \Rd : |z_j| < 1\}$. Recalling the symmetries of $X$ from~\eqref{E : example symmetry}, given $s_0, t_0 \in \R$, positive $\alpha, \beta$, and $\bar y \in \Rk$, we define the paraball
$$
B(s_0, t_0, \bar y, \alpha, \beta) := \phi_{(0, \bar y)} \circ \phi_{(s_0, t_0)} \circ \phi_{(\alpha, \beta)} \big(B(0, 0, 0, 1, 1)\big),
$$
and the corresponding dual paraball 
$$
B^*(s_0, t_0, \bar y, \alpha, \beta) := \psi_{(0, \bar y)} \circ \psi_{(s_0, t_0)} \circ \psi_{(\alpha, \beta)} \big(B^*(0, 0, 0, 1, 1)\big).
$$
\end{definition}

Put differently, setting $\bar x := \bar y + s_0 \gamma_0(t_0)$, the paraball $B(s_0, t_0, \bar y, \alpha, \beta)$ is the set of points $(s, x) \in \R \times \Rk$ satisfying
\begin{itemize}
\item$ |s - s_0| < \alpha$;
\item $|\sum_{i = 1}^m \binom{m}{i}{(- t_0)}^{m - i}(x_i - \bar x_i) + (- t_0)^m(s - s_0)| < \alpha \beta^m, \quad m \leq d - 1$,     
\end{itemize}
whereas the dual paraball $B^*(s_0, t_0, \bar y, \alpha, \beta)$ consists of $(t, y) \in \R \times \Rk$ satisfying
\begin{itemize}
\item $|t - t_0| < \beta$;
\item $|\sum_{i = 1}^m \binom{m}{i}{(- t_0)}^{m - i}(y_i - \bar y_i) + s_0 (t - t_0)^m| < \alpha \beta^m, \quad m \leq d - 1$.      
\end{itemize}

These paraballs are analogous to the Carnot-Carath\'eodory balls in \cite{NagelSteinWanger1985} and two parameter Carnot-Carath\'eodory balls of~\cite{Tao-Wirght2001}. These natural collection play a key role in obtaining $L^p$ improving estimates (e.g.~\cite{CNSW}, \cite[Lemma $4$]{Gressman2019}, \cite[Section $6.2$]{Gressman2024}, \cite[Lemma $6.1$]{Seeger1998}, \cite[Section $10$]{Plamen-Uhlmann2012}, \cite[Proposition $1.5$]{Tao-Wirght2001}), and obtaining precompactness of extremizing sequences for various averaging operators (e.g.~\cite{Christ_extremal, Christ_quasiextremal, Drouot2014}).

We note that as $\|\rchi_{B^*}\|_{\theta'} = \beta^{1/q_\theta'} (\alpha \beta^{d/2})^{(d - 1)/r_\theta'}$, by symmetry of $X$,
\begin{equation}\label{E : quasi paraball}
\scriptX(B, B^*) \gtrsim |B|^{\frac 1{p_{\theta}}} \|\rchi_{B^*}\|_{\theta'}.
\end{equation}
In other words, each pair $(B, B^*)$ is a $c_\theta$-quasiextremal pair. The following theorem justifies the term ``natural" quasiextremal pair for the paraballs.
\begin{theorem}\label{T : reduce to paraball}
Let $\theta \in (0, 1)$ and $0 < \eps < 1/2$. For each $\eps$-quasiextremal pair $(E, F)$ there exists a paraball $B$ with $|B| \leq |E|$ and $\|\rchi_{B^*}\|_{\theta'} \leq \|\rchi_{F}\|_{\theta'}$ such that
$$
\scriptX(E \cap B, F \cap B^*) \gtrsim \eps^C \scriptX(E, F).
$$
\end{theorem}
\begin{proof}
At the heart of the proof is an application of the method of refinement of Christ. But to be able to apply this, we need substantial arguments to deal with the presence of a mixed norm. This is achieved through a combination of a partition argument from~\cite{BS4} and the Lorentz space estimate Theorem~\ref{T : L(q, s) lorentz}.

As we try to approximate $E$ by paraballs, at various stages we will modify $E$ by replacing it with a large subset, which for notational convenience will continue to be denoted by $E$. The following is the first of these. Without loss of generality, we may assume that 
$$
\beta |E| := \scriptX(E, F) = \eps |E|^{1/p_{\theta}} \|\rchi_F\|_{\theta'}.
$$
For the subset $E_1 := \{(s, x) \in E : X^* (\rchi_F) (s, x) > \lambda \eps^{- p_\theta} \beta\}$, we have $\lambda \eps^{- p_\theta} \beta|E_1| < \scriptX(E_1, F ) \leq \beta |E|$. Thus, $|E_1| \leq |E| \eps^{p_\theta}/\lambda$ and so for $\lambda$ large enough
$$
\scriptX(E_1, F) \lesssim |E_1|^{1/p_\theta} \|\rchi_F\|_{\theta'} \leq |E|^{1/p_\theta} \|\rchi_F\|_{\theta'} \eps/2 \leq \scriptX(E, F)/2.
$$
Since we can afford to lose half of the incidences, we can replace $E$ by $E \setminus E_1$ and assume that $X^* (\rchi_F)\lesssim \eps^{- p_\theta} \beta$ on $E$. Likewise, we may assume that $\beta \lesssim X^* (\rchi_F) \lesssim \eps^{- p_\theta} \beta$ on $E$. In other words, for $T_{(s, x)} := \{t \in \R : x - s \gamma_0 (t) \in F\}$, it holds that
\begin{equation}\label{E : length of Ts}
\beta \lesssim |T_{(s, x)}| \lesssim \eps^{-p_\theta} \beta, \qquad (s, x) \in E.
\end{equation}
We claim that each $T_{(s, x)}$ can be approximated by an interval of comparable length up to a factor of $\eps^{- C}$.

To achieve this interval approximation, we use a partition argument from~\cite[Proposition $5.1$]{BS4}. Our first step is to apply~\cite[Lemma $5.4$]{BS4} : fixing $\delta \in (0, 1)$, given the measurable set $T_{(s,x)}$ there exists an interval $I_{(s, x)}$ of length $\eps^{- m_{(s, x)}} \beta$ with $m_{(s, x)} \in \Z$ so that the following holds. For any subinterval $I' \subset I_{(s, x)}$ of length $\eps^{- m_{(s, x)}} \beta/2$, the measure of the complement $|(I_{(s, x)} \setminus I') \cap T_{(s, x)}| \geq \eps^{\delta m_{(s, x)}} \beta$. We will show that the interval $I_{(s, x)}$ meets the above criteria.

We now decompose $E$ into $E^m := \{(s, x) \in E : m \leq m_{(s, x)}< m + 1\}, m \in \Z$, and construct further subsets $E_j^m : = \{(s, x) \in E^m : I_{(s, x)} \subset J_j^m\}$ where $J_j^m := \bigcup_{k = j - 1}^{j + 1} [k \eps^{- m} \beta, (k + 1) \eps^{- m} \beta]$. Finally, let $F_j^m := \{(t, y) \in F : t \in J_j^m\}$. We now recall the following result relating the number of incidences of each pair in this partition to that of $E, F$.

\begin{lemma}{\bf(\cite[Lemma $5.5$]{BS4})}\label{L : lb}
Let $\delta < 1$ be positive. There exists a fixed positive constant $C_\delta$ such that for each $\eps$-quasiextremal pair $(E, F)$
$$
\scriptX(E_j^m, F_j^m) \leq C_\delta \eps^{m (\frac 1{p_\theta} - \frac 1{r_\theta}) \big(1 + \frac{(d + 2) (d - 1)}2 (\frac 1{\theta} - 1) - \delta (d - 2)\big)} |E_j^m|^{1/p_\theta} |F_j^m|^{1/r_\theta'} |J_j^m|^{1/q_\theta' - 1/r_\theta'}.
$$
\end{lemma}

We first consider $\theta_0 \leq \theta < 1$. In this range, $r_\theta' \leq q_\theta'$ and so an application of H\"older's inequality gives $\|\rchi_{F_j^m}\|_{\theta'} \geq |F_j^m|^{1/r_\theta'} |J_j^m|^{1/q_\theta' - 1/r_\theta'}$. We choose $\delta = \delta_\theta := \frac{(d + 2) (d - 1)}{2 (d -2)} (\frac 1{\theta} - 1)$. Thus, for each fixed $m \in \Z$, by H\"older's inequality, bounded overlap of the projection of $F_j^m$ onto the $t$-coordinate, and $p_\theta < q_\theta$,
\begin{align*}
\scriptX(E^m, F) & \lesssim \eps^{m (\frac 1{p_\theta} - \frac 1{r_\theta})} (\sum |E_j^m|)^{1/{p_\theta}} (\sum \|\rchi_{F_j^m}\|_{\theta'}^{p_\theta'})^{1/{p_\theta'}} \lesssim \eps^{m (\frac 1{p_\theta} - \frac 1{r_\theta})} |E^m|^{\frac 1{p_\theta}} \|\rchi_{F}\|_{\theta'}.
\end{align*}
Recalling that $p_\theta < r_\theta$, let $M$ be such that $\eps^{M (1/{p_\theta} - 1/{r_\theta})} < \eps/2$. Summing over the integers $m \geq M$, and then use H\"older's inequality once more we get
\begin{align}\label{E : small tail for E, F}
\sum_{m \geq M} \scriptX(E^m, F) \lesssim \eps^{M (1/{p_\theta} - 1/{r_\theta})} |E|^{1/{p_\theta}} \|\rchi_{F}\|_{\theta'} < \scriptX(E, F)/2.
\end{align}
In conclusion, for $\theta_0 \leq \theta < 1$, replacing $E$ by $E^C$, we arrive at our desired conclusion that for each $(s, x) \in E$ the interval $I_{(s, x)}$ has length at most $\lesssim \eps^{- C} \beta$.

Now we address the case of $0 < \theta < \theta_0$. As $q_\theta' < r_\theta'$ in this range, the inequality~\eqref{E : small tail for E, F} does not follow directly from Lemma~\ref{L : lb}. To overcome this, using the Lorentz space estimate we will first do a bit of pruning of the set $F$. With $F^k := F \cap \pi^{- 1}(\{t : 2^k \leq  \int \rchi_F (t, y) dy < 2^{k + 1}\}), k \in \Z$, Theorem~\ref{T : L(q, s) lorentz} implies the existence of one $F^k$ so that
$$
\|\rchi_{F^k}\|_{\theta'} \gtrsim \eps^C \|\rchi_F\|_{\theta'} \qtq{and} \scriptX(E, F^k) \gtrsim \eps^C \scriptX(E, F).
$$

Since we can afford to lose a factor of $\eps^C$, we can replace $F$ by one such $F^k$ (and as before continue to denote it by $F$). By construction, $\|\rchi_F\|_{\theta'} \sim |F|^{1/r_\theta'} |\pi(F)|^{1/q_\theta' - 1/r_\theta'}$. Thus, by Lemma~\ref{L : lb}, with $B_\delta := \big(1 + \frac{(d + 2)(d - 1)}{2}(\frac{1}{\theta} - 1) - \delta (d - 2)\big)(\frac1{p_\theta} - \frac1{r_\theta})$,
$$
\scriptX(E_j^m, F_j^m) \lesssim \eps^{m (B_\delta - 1/q_\theta' + 1/r_\theta')} |E_j^m|^{\frac{1}{p_\theta}} \|\rchi_{F_j^m}\|_{\theta'}.
$$
For $0 < \theta < \theta_0$, we note that at $\delta = 0$ the exponent $B_\delta - 1/q_\theta' + 1/r_\theta'$ is positive. This lets us choose a positive $\delta = \delta(\theta)$ such that $B_\delta - 1/q_\theta' + 1/r_\theta' > 0$ and so
$$
\scriptX(E_j^m, F_j^m) \lesssim \eps^{c m} |E_j^m|^{1/p_\theta} \|\rchi_{F_j^m}\|_{\theta'}.
$$
We now proceed as in $\theta_0 \leq \theta < 1$ and repeat the argument leading to~\eqref{E : small tail for E, F}. Thus, our goal of showing that on $E$ each length $|I_{(s, x)}| \lesssim \eps^{- C} \beta$ is finally achieved.

Next, recalling the sets $E_j (:= E_j^C)$ and $F_j (:= F_j^C)$, we may assume that there exists an $\sim \eps$-quasiextremal pair $(E_j, F_j)$ with $|E_j| \gtrsim \eps^C |E|$ and $\|\rchi_{F_j}\|_{\theta'} \gtrsim \eps^C \|\rchi_F\|_{\theta'}$. For otherwise, with $J_\mu := \{j : \scriptX(E_j, F_j) < \mu \eps |E_j|^{1/p_\theta} \|\rchi_{F_j}\|_{\theta'}\}$, by H\"older's inequality (using that $p_\theta < q_\theta$),
\begin{align*}
\sum_{j \in J_\mu} \scriptX(E_j, F) \leq \mu \eps (\sum_j |E_j|)^{1/p_\theta} (\sum_j \|\rchi_{F_j}\|_{\theta')}^{q_\theta'})^{1/q_\theta'} \lesssim \mu \scriptX(E, F).
\end{align*}
Choosing a positive $\mu$ small enough gives an $\sim \eps$-quasiextremal pair $(E_j, F_j)$. Furthermore, for the favorable collection of $E_j$, once more by H\"older's inequality,
\begin{align*}
\eps |E|^{\frac 1{p_{\theta}}} \|\rchi_F\|_{\theta'} = \scriptX(E, F) \lesssim (\sum |E_j|^{\frac{q_\theta}{p_\theta}})^{\frac 1{q_\theta}} (\sum \|\rchi_{F_j}\|_{\theta'}^{q_\theta'})^{\frac 1{q_\theta'}} \lesssim |E|^{\frac 1{q_\theta}} \|\rchi_F\|_{\theta'} \sup_j |E_j|^{\frac 1{p_\theta} - \frac 1{q_\theta}}.
\end{align*} 
This guarantees the existence of one $E_j$ satisfying $|E_j| \gtrsim \eps^C |E|$. Similarly, for $F_j$. In summary, replacing $E$ by one such $E_j$ (and $F \rightsquigarrow F_j$) we may conclude that there exists a fixed interval $J$ of length $\sim \eps^{- C} \beta$ (we recall that the sets $\pi(F), T_{(s, x)}$ are of measure $\gtrsim \beta$) so that
\begin{equation}\label{E : find interval J}
\pi(F) \subset J \qtq{and} I_{(s, x)} \subset J \qtq{for each} (s,x) \in E.
\end{equation}

We are at the final stage of the proof. We give the details for $d = 2 D + 1$ odd, the other case being identical. Applying the method of refinement using property~\eqref{E : find interval J}, we get $(t_0, \bar y) \in F$ and $(s_0, \bar x)\in E$ with $\bar x = \bar y + s_0 \gamma_0(t_0)$, so that for $\alpha|F| := \scriptX(E, F)$, the following holds. There exist measurable sets $\Omega_1, \Omega_2 \subset \Rd$ with $|\Omega_1| \sim \alpha^{D + 1} \beta^D$ and $|\Omega_2| \sim \alpha^D \beta^{D + 1}$, such that for each $(s, x) := (s_1, t_1, \ldots, s_D, t_D, s_{D + 1}) \in \Omega_1$
\begin{itemize}
\item $\Phi (s, x) := (s_{D + 1}, \bar y + s_1 \gamma_0(t_0) - \sum_{j = 1}^D (s_j - s_{j - 1}) \gamma_0(t_j)) \in E$;
\item $|s_i - s_{i - 1}| \lesssim \eps^{- C} \alpha$, \, $|t_j-t_i|\lesssim \eps^{- C} \beta$,
\end{itemize}
and for each $(t, y) := (t_1, s_1, \ldots, t_D, s_D, t_{D + 1}) \in \Omega_2$
\begin{itemize}
\item $\Psi (t, y) := (t_{D + 1}, \bar x - \sum_{j = 1}^D (s_{j - 1} - s_j) \gamma_0(t_j) - s_D \gamma_0(t_{D + 1})) \in F$;
\item $|s_i - s_{i - 1}| \lesssim \eps^{- C} \alpha$, \, $|t_j - t_i| \lesssim \eps^{- C} \beta$.
\end{itemize}

Fixing a large $C$, we consider the paraball $B = B(s_0, t_0, \bar y, C \eps^{- C} \alpha, C \eps^{- C} \beta)$ and its dual. We claim $\Phi(\Omega_1) \subset B$. To prove this, by symmetry, we may assume that $\eps = 1, s_0 = t_0 = 0$ and $\bar y = \bar x = 0 \in \R^{d - 1}$. Thus,
$|s_{D + 1}| \lesssim \eps^{- C} \alpha$ and
$$
|\sum_{1 \leq j \leq D} (s_j - s_{j - 1}) t_j^m| \lesssim \eps^{- C} \alpha \beta^m
$$
for each $m \leq d - 1$, establishing our claim. Similarly, $\Psi(\Omega_2) \subset B^*$. An application of the following Lemma~\ref{L : partition by parablls} finally finishes the proof of Theorem~\ref{T : reduce to paraball}.
\end{proof}

\begin{lemma}\label{L : partition by parablls}
Given a paraball $B$ and a positive $\delta < 1$, there exists a family of paraballs $\{B_l : l \in L\}$ with $|L| \lesssim \delta^{- C}$ so that $B \subset \bigcup_l B_l$ and $B^{*} \subset \bigcup_l B_l^{*}$ with
$$
|B_l| \sim \delta |B| \qtq{and} \|\rchi_{B^*_l}\|_{\theta'} \sim \delta \|\rchi_{B^*}\|_{\theta'}.
$$
\end{lemma}
\begin{proof}
Applying symmetry, we may assume that $B$ is the unit cube in $\Rd$. Fixing $\eta_1 > 0$ to be chosen precisely for later purpose, let $\eta_1^d \eta_2^{(d^2 - d)/2} := \delta$. Let $\{\bar y^i : i \in I\}$ denote a maximal $\eta_1 \eta_2^d$-separated subset of the unit cube in $\Rk$. Next, we choose maximal $\eta_1$ and $\eta_2$-separated subsets of $[- 1, 1]$ that we denote by $\{s^j : j \in J\}$ with $s^1 := 0$ and $\{t^k : k \in K\}$ with $t^1 := 0$, respectively. Thus, $|I| + |J| + |K| \lesssim (\eta_1 \eta_2)^{- C}$.

For the paraballs $B_{i, j, k} := B(s^j, t^k, \bar y^i, 2\eta_1, 2\eta_2)$, clearly $|B_{i, j, k}| \sim \eta_1^d \eta_2^{(d^2 - d)/2} = \delta$ and $\|\rchi_{B^{*}_{i, j, k}}\|_{\theta'} \sim \eta_2^{1/q_\theta'} (\eta_1^{d - 1} \eta_2^{(d^2 - d)/2})^{1/r_\theta'}$. Now we specify $\eta_1$ so that $\|\rchi_{B^{*}_{i, j ,k}}\|_{\theta'} \sim \delta$ by declaring $\eta_2^{\frac 1{{q_\theta'}} + \frac{d - 1}{2 r_\theta'}} := \delta^{\frac 1{r_\theta} + \frac 1{d r'_\theta}}$. From Definition~\ref{E : def paraball} it directly follows that
\begin{align*}
& B \subset \bigcup_{i, k} B(0, t^k, \bar y^i, 2\eta_1, 2\eta_2) \subset \bigcup_{i, j, k} B_{i, j, k}
\\
& B^* \subset \bigcup_{i, j} B^*(s^j, 0, \bar y^i, 2 \eta_1, 2 \eta_2) \subset \bigcup_{i, j, k} B_{i, j, k}^*.
\end{align*}
This completes the proof.
\end{proof}

Theorem~\ref{T : reduce to paraball} in conjunction with Theorem~\ref{T : L(p, s) Lorentz} and Theorem~\ref{T : L(q, s) lorentz} implies the following immediate corollary.

\begin{corollary}\label{C : single ball}
Let $\theta \in (0, 1)$ and $\eps$ be positive. For each $f = \sum_j 2^j f_j, f_j \sim \rchi_{E_j}$ satisfying $\|Xf\|_{L_t^{q_\theta} L_y^{r_\theta}} \geq \eps \|f\|_{L^{p_\theta}}$, there exist an index $j_0 \in \Z$ and a paraball $B$ so that
$$
2^{j_0} |E_{j_0} \cap  B|^{1/p_\theta} \gtrsim \eps^C \|f\|_{L^{p_\theta}} \qtq{and} |B| \leq |E_{j_0}|.
$$
\end{corollary}

\section{Mock-distance}\label{S : mock distance}

In order to establish the uniform localization of a near extremizer, Proposition~\ref{P : loc}, we need to analyze the interaction between the image of two distant paraballs under $X$. To quantify this interaction, we consider a notion of distance on the set of paraballs, in the spirit of the pseudo-distance in~\cite{Christ_extremal}.

For two paraballs $B^a = B(s_a, t_a, \bar y^a, \alpha_a, \beta_a), B^b = B(s_b, t_b, \bar y^b, \alpha_b, \beta_b)$ (recalling our convention $\bar x := \bar y + s \gamma_0(t)$) let $\rho$ be defined by
\begin{align}\label{E : def dist}
\rho (B^a, B^b) &\notag : = (\frac{\alpha_a}{\alpha_b} + \frac{\alpha_b}{\alpha_a}) + (\frac{\beta_a}{\beta_b} + \frac{\beta_b}{\beta_a}) + |s_a - s_b| (\frac 1{\alpha_a} + \frac 1{\alpha_b}) + |t_a - t_b| (\frac 1{\beta_a} + \frac 1{\beta_b})
\\
&\notag + \sum_{m = 1}^{d - 1} (\alpha_a \beta_a^m)^{- 1} \Big|\sum_{i = 1}^m \binom{m}{i}{(- t_a)}^{m - i}(\bar x^b_i - \bar x^a_i) + (- t_a)^m (s_b - s_a)\Big|
\\
&\notag + \sum_{m = 1}^{d - 1} (\alpha_b \beta_b^m)^{- 1} \Big|\sum_{i = 1}^m \binom{m}{i}{(- t_b)}^{m - i}(\bar x^a_i - \bar x^b_i) + (- t_b)^m (s_a - s_b)\Big|
\\
&\notag + \sum_{m = 1}^{d - 1} (\alpha_a \beta_a^m)^{- 1} \Big|\sum_{i = 1}^m \binom{m}{i}{(- t_a)}^{m - i}(\bar y^b_i - \bar y^a_i) + s_a (t_b - t_a)^m\Big|
\\
& + \sum_{m =1 }^{d - 1} (\alpha_b \beta_b^m)^{- 1} \Big|\sum_{i = 1}^m \binom{m}{i}{(- t_b)}^{m - i}(\bar y^a_i - \bar y^b_i) + s_b (t_a - t_b)^m\Big|.
\end{align}

A few comments are in order. $\rho$ is not a distance function as $\rho (B, B) = 4$. But this is of little significance to us, as we will primarily be interested in the behavior of the image of two distant (in $\rho$) paraballs under $X$. As it is not a pseudo-distance either, for lack of a better term, we call it a mock-distance.

Next, we take a close look at the terms defining $\rho$. The first (second) term compares the length of the base of the paraballs (dual paraballs). The third (fourth) term measures the distance between the first coordinate of the center of the paraballs (dual paraballs). The fifth and sixth terms measure the distance between the last $(d - 1)$ coordinates of the center of the paraballs, and likewise the seventh and eighth terms corresponding to the dual paraballs. We will observe momentarily that the seventh and eighth terms are redundant as these are essentially dominated by the fifth and sixth terms, respectively. However, we still include these redundant terms to ensure that the mock-distance is dually symmetric, i.e. $\rho(B^a, B^b) = \rho({B^a}^*, {B^b}^*)$ which will be useful in our analysis.

We note down the following immediate properties of $\rho$ whose routine proof we omit (see, e.g., ~\cite[section $5$]{CB}).
\begin{itemize}
\item $\rho$ is preserved under the symmetry $\phi$ of $X$ : $\rho (B_1, B_2) = \rho (\phi(B_1), \phi(B_2))$.
\item The following almost-triangle inequality holds.
\begin{equation}\label{E : almost trian}
\rho (B_1, B_2) \lesssim \rho (B_1, B_3)^C + \rho (B_2,B_3)^C;
\end{equation}
\item The following covering property holds.
\begin{equation}\label{E : scaling covers}
B_1 \subset C \rho (B_1, B_2)^C B_2
\end{equation}
where $\lambda B$ denotes the paraball with same center but the width parameters $\alpha, \beta$ scaled by $\lambda$. Note that this does not correspond to a symmetry of $X$.
\end{itemize}
Furthermore, $\rho$ satisfies~\eqref{E : dist and volume}, which roughly says that two $\rho$-distant paraballs of comparable volume have a small intersection regardless of their width parameters.

\begin{proposition}\label{P : dist and intersection}
There exists a fixed positive constant $c$ such that
\begin{equation}\label{E : dist and volume}
|B^a \cap B^b| \lesssim \rho (B^a, B^b)^{- c} \max(|B^a|, |B^b|).
\end{equation}
\end{proposition}

\begin{proof}
The proof follows an argument in~\cite{Christ_extremal}. Suppose that $B^a = B(s_a, t_a, \bar y^a, \alpha_a, \beta_a)$ and $B^b = B(s_b, t_b, \bar y^b, \alpha_b, \beta_b)$ with $\alpha_a \leq \alpha_b$. In case the first term in~\eqref{E : def dist} dominates $\rho$, with $S := (s_a - \alpha_a, s_a + \alpha_a) \cap (s_b - \alpha_b, s_b + \alpha_b)$,
\begin{align*}
|B^a & \cap B^b| \leq |S| \alpha_b^{- 1} \max(|B^a|, |B^b|) \leq \alpha_a \alpha_b^{- 1} \max(|B^a|, |B^b|).
\end{align*}
Thus, bound~\eqref{E : dist and volume} follows. Similarly, if the second term dominates $\rho(B^a, B^b)$ and the first term $\frac{\alpha_a}{\alpha_b} + \frac{\alpha_b}{\alpha_a} \leq c \rho (B^a, B^b)^c$, with $A := \max(\alpha_a^2 \beta_a^d, \alpha_b^2 \beta_b^d)/\min(\alpha_a^2 \beta_a^d, \alpha_b^2 \beta_b^d)$,
\begin{align*}
(\beta_a/\beta_b + \beta_b/\beta_a)^d \lesssim \max(\alpha_b^2 \beta_a^d, \alpha_b^2 \beta_b^d\big)/\min\big(\alpha_a^2 \beta_a^d, \alpha_a^2 \beta_b^d) \leq A (c' \rho (B^a, B^b))^C.
\end{align*}
Thus, if $\alpha_a^2 \beta_a^d \leq \alpha_b^2 \beta_b^d$, the required upper bound follows from
$$
|B^a \cap B^b| \leq (\alpha_a^2 \beta_a^d)^{(d - 1)/2} \min(\alpha_a, \alpha_b) \leq A^{(1 - d)/2} \max(|B^a|, |B^b|).
$$

For the interval $S$ as above, $|S| \lesssim (|s_a - s_b| \big(1/\alpha_a + 1/\alpha_b\big))^{- 1} \max(\alpha_a, \alpha_b)$, so if the third term in~\eqref{E : def dist} dominates $\rho$, bound~\eqref{E : dist and volume} follows from
$$
|B^a \cap B^b| \lesssim \Big[|s_a - s_b| \big(1/\alpha_a + 1/\alpha_b\big)\Big]^{- 1} \max(|B^a|, |B^b|).
$$

Now we consider the fourth term in~\eqref{E : def dist}. For $(s, x) \in B^a \cap B^b$, we have $|x_1 - \bar x^a_1 - t_a (s - s_a)| \leq \alpha_a \beta_a$ and $|x_1 - \bar x^b_1 - t_b (s - s_b)| \leq \alpha_b \beta_b$. With $e := \bar x_1^a - \bar x_1^b - s_a t_a + s_b t_b$, subtracting gives $|s (t_a - t_b) + e| \leq 2 \max (\alpha_a \beta_a, \alpha_b \beta_b)$. Thus, if $|t_a - t_b| \gtrsim \rho (B^a, B^b) \beta_a$, the measure of the set of all such $s$ is $\lesssim \rho (B^a, B^b)^{- 1} \alpha_a$, producing the desired upper bound on $|B^a \cap B^b|$.

Next, let us consider when the fifth term in~\eqref{E : def dist} dominates, thus, for some $m$,
$$
\Big|\sum_{i = 1}^m \binom{m}{i} {(- t_a)}^{m - i}(\bar x^b_i - \bar x^a_i) + (- t_a)^m (s_b - s_a)\Big| \gtrsim \alpha_a \beta_a^m \rho (B^a, B^b).
$$ 
For $j \leq d - 1$, we consider the degree one polynomial $Q_j^b$ on $\R \times \Rk$ given by
$$
Q_j^b(s, x) := \sum_{i = 1}^j \binom{j}{i}{(- t_b)}^{j - i} (x_i - \bar x_i^b) + (- t_b)^j (s - s_b).
$$
Given $s \in \R$, we consider $x (s) \in \Rk$ so that $Q_j^b(s, x(s)) = 0$ for each $j \leq d - 1$ (using the convention $0^0 := 1$). For the degree one polynomial $P(s) := Q_m^a (s, x(s))$, where $Q_m^a$ is defined by replacing $b$ with $a$ in the terms defining $Q_m^b$, we note that 
$$
|P(s)| > \alpha_a \beta_a^m + \alpha_b \beta_b^m \qtq{implies that} B^a \cap B^b\cap 
(\{s\} \times \R^{d - 1}) = \emptyset.
$$
Let $\eps_0 := \frac{2 \max (\alpha_a \beta_a^m, \alpha_b \beta_b^m)}{\rho (B^a, B^b) \alpha_a \beta_a^m}$. Since the first two terms in~\eqref{E : def dist} do not dominate $\rho^{\frac 1{2 d}}$, we have $\eps_0 \lesssim \rho (B^a, B^b)^{- 1/3}$. Our assumption is that $|P(s_b)| \gtrsim \alpha_a \beta_a^m \rho (B^a, B^b)$. Thus, for $s \in (s_b - \alpha_b, s_b + \alpha_b)$ except on a subset of measure $\lesssim \eps_0 \alpha_b$,
$$
|P(s)| \geq \eps_0 \rho (B^a, B^b) \alpha_a \beta_a^m = 2 \max(\alpha_a \beta_a^m, \alpha_b \beta_b^m) \geq \alpha_a \beta_a^m + \alpha_b \beta_b^m.
$$
In consequence, we arrive at the bound $|B^a \cap B^b| \lesssim \eps_0 |B^b| \lesssim \rho (B^a, B^b)^{- 1/3} |B^b|$. Identical arguments give the bound~\eqref{E : dist and volume} when the sixth term dominates $\rho$.

Looking at the the seventh term in~\eqref{E : def dist}, suppose that for some $m \leq d - 1$, 
$$
(\alpha_a \beta_a^m)^{- 1} \Big|\sum_{i = 1}^m \binom{m}{i}{(- t_a)}^{m - i}(\bar y^b_i - \bar y^a_i) + s_a (t_b - t_a)^m\Big| \gtrsim \rho (B^a, B^b),
$$
and each of the preceding six terms is $\leq c' \rho (B^a, B^b)$ where $c' > 0$ is sufficiently small to be chosen momentarily. We will show that this implies that the sixth term also dominates $\rho$ and thus produces a contradiction. By symmetry, we may assume that $B^a = B(0, 0, 0, 1, 1)$ and $B^b = B(s_b, t_b, \bar y^b, \alpha_b, \beta_b)$. Our hypotheses is that
\begin{equation}\label{E : large eigth term}
|\bar y^b_m| > c\rho (B^a, B^b),
\end{equation}
and additionally, $|s_b| + |\bar x^b_m| < c' \rho (B^a, B^b), |t_b| < c' \rho (B^a, B^b)^{1/m}$. But then $|\bar y^b_m| = |\bar x^m_b - s_b t^m_b| < C c' \rho (B^a, B^b)$ and so choosing $c'$ small enough contradicts~\eqref{E : large eigth term}. Identical arguments produce~\eqref{E : dist and volume} when the eighth and final term of~\eqref{E : def dist} dominates the mock distance $\rho$. This completes the proof.
\end{proof}

\section{Proof of Proposition~\ref{P : loc} : Localization of $X$}\label{S : proof of localization}

We now have the tools required to prove Proposition~\ref{P : loc}.

\begin{proof}[Proof of Proposition~\ref{P : loc}] Assuming the existence of a symmetry $\phi$ we deduce its post hoc independence from $\eps$. We may assume that $f$ is normalized and $\|X f\|_\theta \geq A_\theta/2$. By Corollary ~\ref{C : single ball} there exists a paraball $B$ and an index $j \in \Z$ so that
$$
\|2^j f_j \rchi_B\|_{L^{p_\theta}} \sim 1.
$$
By symmetry, we may assume that $j = 0$ and the paraball $B$ is the unit cube in $\Rd$. Now suppose that for a symmetry $\phi = \phi_{(0, \bar y)} \circ \phi_{(s_0, t_0)} \circ \phi_{(\alpha, \beta)}$ it holds that
$$
\|\phi^* f\|_{L^{p_\theta}(\{|(s, x)| > R\} \cup \{|\phi^*f| > R\})} < \eps
$$
for a sufficiently small $\eps$. We will prove that the above also holds for $\phi$ equals the identity map on $\Rd$ with $R$ replaced by a possibly larger $\tilde R$.

By our assumption $\|\phi^* f_0\|_{L^{p_\theta}(\{|(s, x)| > R\} \cup \{|\phi^*f_0| > R\})} < \eps$. For $\eps$ sufficiently small, by the triangle inequality and unwinding the notation, we have
$$
\|f_0\|_{L^{p_\theta}(B \cap \{|f| \lesssim 1\} \cap \{|f_0| < (\alpha \beta^{(d - 1)/2})^{- d/{p_\theta}} R\}) \cap F_1 \cap F_2)} \sim 1
$$
where $F_1 := \{(s, x) : |s - s_0| < 2 \alpha R\}$ and $F_2 := \{(s, x) :|\big(G_{t_0}^{- 1} (x - (\alpha s + s_0) \gamma_0(t_0) - \bar y)\big)_i| < 2 \alpha \beta^i R, \, \text{for} \, 1 \leq i \leq d - 1\}$. By repeated application of H\"older's inequality,
\begin{align}\label{E : aplha R > 1}
&  1 \lesssim \|f_0\|_{L^{p_\theta}(\{|f| \lesssim 1\} \cap F_1 \cap B)} \lesssim \alpha R,
\\\label{E : alpha beta R > 1}
& 1 \lesssim \|f_0\|_{L^{p_\theta}(\{|f| \lesssim 1\} \cap F_2 \cap B)} \lesssim \alpha \beta^i R, \quad 1 \leq i \leq d - 1,
\\\label{E : alphe^-1 R > 1}
& 1 \lesssim \|f_0\|_{L^{p_\theta}(B \cap \{|\phi^*f_0| < R\})} \lesssim (\alpha \beta^{(d - 1)/2})^{- d/{p_\theta}} R.
\end{align}

For $\beta \geq 1$, combining ~\eqref{E : aplha R > 1} and\eqref{E : alphe^-1 R > 1} and a bit of arithmetic, we get
$$
R^{- 1} \lesssim \alpha \lesssim R^{\frac{p_\theta}d} \qtq{and so} \beta \sim 1.
$$
Similarly, considering $B, \{|f| \lesssim 1\}, F_2$, by H\"older's we get $|s_0| + |t_0| + |\bar y| \lesssim R^{1 + \frac{p_\theta}d}$. In consequence, 
$$
\{\phi(s, x) : |(s, x)| > R\} \supset \{|(s, x)| \gtrsim R^{d - 1 + p_\theta}\} \qtq{and} \{|\phi^*f| > R\} \supset \{|f| \gtrsim R^{1 + \frac d{p_\theta}}\}.
$$

Identical arguments work for $\beta \leq 1$. Indeed, combining~\eqref{E : alpha beta R > 1} at $i = d - 1$ and~\eqref{E : alphe^-1 R > 1} gives
$$
1 \lesssim \alpha \beta^{d  -1} R \lesssim R^{\frac {p_\theta}d + 1} \beta^{(d - 1)/2}.
$$
Thus, $R^{- (1 + \frac {p_\theta}d) \frac 2{d - 1}} \lesssim \beta \leq 1$. Now we proceed exactly as above and apply H\"older's inequality a few times to get 
$$
R^{- 1} \lesssim \alpha \lesssim R^{1 + \frac {2 p_\theta}d} \qtq{and} |s_0| + |t_0| + |\bar y| \lesssim R^{2 + \frac{2 p_\theta}d}.
$$
For $\beta \leq 1$ this yields
$$
\{\phi(s, x) : |(s, x)| > R\} \supset \{|(s, x)| \gtrsim R^{2(d - 1 + p_\theta)}\}, \, \text{and} \, \{|\phi^*f| > R\} \supset \{|f| \gtrsim R^{2 (1 + \frac d{p_\theta})}\}.
$$
Finally, combining the above two inclusions, with $\tilde R :=  C R^{2 (d + p_\theta + \frac d{p_\theta})}$, we obtain 
$$
\|f\|_{L^{p_\theta}(\{|(s, x)| > \tilde R\} \cup \{|f| > \tilde R\})} < \eps.
$$

We now return to the main thread of the proof. It uses the method of~\cite{Christ_extremal}. We sketch the argument below; complete and identical details are given in~\cite[Lemma~$11.2$]{CB}. Let $\eps > 0$ be fixed for the remainder of the proof. By Proposition~\ref{P : freq loc} there exists positive $\delta$ such that the following holds for each normalized $f$. If $f = \sum 2^j f_j, f_j \sim \rchi_{E_j}$ for a mutually disjoint collection $\{E_j\}$ and satisfies $\|Xf\|_{L_t^{q_\theta} L_y^{r_\theta}} \geq (1 - \delta) A_\theta$, then there exists a collection of frequencies $S \subset \Z$ contained in an interval of length $\sim \eps^{- C}$ so that for the truncated function $\|X (\sum_{j \in S} 2^j f_j)\|_{L_t^{q_\theta} L_y^{r_\theta}} \geq (1 - \eps) A_\theta$. Next, we apply Corollary~\ref{C : single ball} to this truncated function and get a distinguished frequency $j_0 \in S$ and a paraball $B_{j_0}$ so that 
$$
\|2^{j_0} f_{j_0} \rchi_{B_{j_0}}\|_{L^{p_\theta}} \sim 1 \qtq{and} |B_{j_0}| \leq |E_{j_0}|.
$$
Applying suitable symmetry we may assume that $j_0 = 0$ (and so the frequency support $S$ is contained in an interval in $\Z$ centered at $j = 0$ of length $\sim \eps^{- C}$) and $B_0$ is the unit cube in $\Rd$. Thus, it suffices to show that for each $j \in S$,
$$
E_{j}\subset C\eps^{-C^{C\eps^{-C}}} B_0.
$$

Next, we recall~\cite[Lemma~$11.1$]{CB}. Analogous statement holds for our paraballs, and the proof is exactly the same as the proof in~\cite{CB} with the only change being incorporation of Corollary~\ref{C : single ball} of this article. Fixing a frequency $j \in S$, we thus collect a collection of paraballs $\{B_{i, j}\}_{i \in N_j}$ with $|N_j| \lesssim \eps^{- C}$ so that $E_j \subset \bigcup_{i \in N_j} B_{i, j}$. Thus, by the engulfing property~\eqref{E : scaling covers}, it suffices to show that each paraball $B_{i, j}$ lies within a $\rho$-distance of $\sim \eps^{- C^{C\eps^{-C}}}$ of the unit cube. To show this, we choose a positive parameter $\lambda$ for later purposes. Suppose, if possible, that the collection of paraballs $B_{i, j}$ can be decomposed into two disjoint classes with the property that each pair consisting of one from each class is at least at a distance $\eps^{- \lambda}$ in $\rho$ from each other (we may assume that $\eps^{- \lambda} \geq 5$). Arguing as in the proof of Proposition~\ref{P : freq loc}, decomposing the function $f$ into two pieces corresponding to these two classes, for $\eps$ sufficiently small, produces a pair of paraballs $(B_1, B_2)$ such that 
$$
\eps^C |B_1| \lesssim |B_1 \cap B_2| \leq |B_2| \lesssim \eps^{- C} |B_1| \qtq{and} \rho (B_1, B_2) \gtrsim \eps^{- \lambda}.
$$
(see, e.g., ~\cite[Lemma~$11.2$]{CB}). But this contradicts the conclusion of Proposition~\ref{P : dist and intersection} if $\lambda$ is chosen large enough depending on $\theta$. In summary, it is not possible to construct such a decomposition of the set of all paraballs $B_{i, j}$. Since there are $\lesssim \eps^{- C}$ many paraballs, applying the almost-triangle inequality~\eqref{E : almost trian}, the proof is complete.
\end{proof}



\bibliographystyle{plain}
\bibliography{Xray_moment}

\begin{thebibliography}{10}

\bibitem{BernsteinLoss1997}
Albert Baernstein and Michael Loss.
\newblock Some conjectures about {$L^p$} norms of {$k$}-plane transforms.
\newblock {\em Rend. Sem. Mat. Fis. Milano}, 67:9--26, 1997.

\bibitem{BennettBezFlock2018}
Jonathan Bennett, Neal Bez, Taryn~C. Flock, Susana Guti\'errez, and Marina
  Iliopoulou.
\newblock A sharp {$k$}-plane {S}trichartz inequality for the {S}chr\"odinger
  equation.
\newblock {\em Trans. Amer. Math. Soc.}, 370(8):5617--5633, 2018.

\bibitem{BennettCarberyChristTao2008}
Jonathan Bennett, Anthony Carbery, Michael Christ, and Terence Tao.
\newblock The {B}rascamp-{L}ieb inequalities: finiteness, structure and
  extremals.
\newblock {\em Geom. Funct. Anal.}, 17(5):1343--1415, 2008.

\bibitem{CB}
Chandan Biswas.
\newblock Existence of extremizers for a model convolution operator.
\newblock {\em Journal of Fourier Analysis and Applications}, 25:2653--2689,
  2019.

\bibitem{Bourgain1986}
Jean Bourgain.
\newblock Averages in the plane over convex curves and maximal operators.
\newblock {\em J. Analyse Math.}, 47:69--85, 1986.

\bibitem{Bourgain1991}
Jean Bourgain.
\newblock Besicovitch type maximal operators and applications to {F}ourier
  analysis.
\newblock {\em Geom. Funct. Anal.}, 1(2):147--187, 1991.

\bibitem{Rubinhyperbolic2019}
William~O. Bray and Boris Rubin.
\newblock Radon transforms over lower-dimensional horospheres in real
  hyperbolic space.
\newblock {\em Trans. Amer. Math. Soc.}, 372(2):1091--1112, 2019.

\bibitem{CarberyChristWaingerWatson1989}
Anthony Carbery, Michael Christ, James Vance, Stephen Wainger, and David~K.
  Watson.
\newblock Operators associated to flat plane curves: {$L^p$} estimates via
  dilation methods.
\newblock {\em Duke Math. J.}, 59(3):675--700, 1989.

\bibitem{CarberyRicciWright2003}
Anthony Carbery, Fulvio Ricci, and James Wright.
\newblock Maximal functions and singular integrals associated to polynomial
  mappings of {$\Bbb R^n$}.
\newblock {\em Rev. Mat. Iberoamericana}, 19(1):1--22, 2003.

\bibitem{CarberySeegerWainger1999}
Anthony Carbery, Andreas Seeger, Stephen Wainger, and James Wright.
\newblock Classes of singular integral operators along variable lines.
\newblock {\em J. Geom. Anal.}, 9(4):583--605, 1999.

\bibitem{Christ1984}
Michael Christ.
\newblock Estimates for the {$k$}-plane transform.
\newblock {\em Indiana Univ. Math. J.}, 33(6):891--910, 1984.

\bibitem{Christ1985groups}
Michael Christ.
\newblock Hilbert transforms along curves. {I}. {N}ilpotent groups.
\newblock {\em Ann. of Math. (2)}, 122(3):575--596, 1985.

\bibitem{Christ_1998}
{M}ichael {C}hrist.
\newblock Convolution, curvature, and combinatorics: a case study.
\newblock {\em Internat. Math. Res. Notices}, 1998(19):1033--1048, 1998.

\bibitem{Christ_extremal}
{M}ichael {C}hrist.
\newblock On extremals for a {R}adon-like transform.
\newblock {\em Preprint, arXiv:1106.0728}, 2011.

\bibitem{Christ_quasiextremal}
{M}ichael {C}hrist.
\newblock Quasiextremals for a {R}adon-like transform.
\newblock {\em Preprint, arXiv:1106.0722.}, 2011.

\bibitem{Christ2014}
Michael Christ.
\newblock Extremizers of a {R}adon transform inequality.
\newblock In {\em Advances in analysis: the legacy of {E}lias {M}. {S}tein},
  volume~50 of {\em Princeton Math. Ser.}, pages 84--107. Princeton Univ.
  Press, Princeton, NJ, 2014.

\bibitem{ChristDendrinosStovallStreet2020}
Michael Christ, Spyridon Dendrinos, Betsy Stovall, and Brian Street.
\newblock Endpoint {L}ebesgue estimates for weighted averages on polynomial
  curves.
\newblock {\em Amer. J. Math.}, 142(6):1661--1731, 2020.

\bibitem{ChristErdogan2002}
Michael Christ and M.~Burak {E}rdo\u{g}an.
\newblock Mixed norm estimates for a restricted {X}-ray transform.
\newblock {\em J. Anal. Math.}, 87:187--198, 2002.

\bibitem{ChristErdogan2008}
Michael Christ and M.~Burak {E}rdo\u{g}an.
\newblock Mixed norm estimates for certain generalized {R}adon transforms.
\newblock {\em Trans. Amer. Math. Soc.}, 360(10):5477--5488, 2008.

\bibitem{CNSW}
Michael Christ, Alexander Nagel, Elias~M. Stein, and Stephen Wainger.
\newblock Singular and maximal radon transforms: Analysis and geometry.
\newblock {\em Annals of Mathematics}, 150(2):489--577, 1999.

\bibitem{ChristXue2012}
Michael Christ and Qingying Xue.
\newblock Smoothness of extremizers of a convolution inequality.
\newblock {\em J. Math. Pures Appl. (9)}, 97(2):120--141, 2012.

\bibitem{DendrinosLaghiWright2009}
Spyridon Dendrinos, Norberto Laghi, and James Wright.
\newblock Universal {$L^p$} improving for averages along polynomial curves in
  low dimensions.
\newblock {\em J. Funct. Anal.}, 257(5):1355--1378, 2009.

\bibitem{BS4}
{S}pyridon {D}endrinos and {B}etsy {S}tovall.
\newblock Uniform estimates for the {X}-ray transform restricted to polynomial
  curves.
\newblock {\em J. Funct. Anal.}, 262:4986--5020, 2012.

\bibitem{Drouot2014}
Alexis Drouot.
\newblock Sharp constant for a {$k$}-plane transform inequality.
\newblock {\em Anal. PDE}, 7(6):1237--1252, 2014.

\bibitem{Druot2015}
Alexis Drouot.
\newblock Quantitative form of certain {$k$}-plane transform inequalities.
\newblock {\em J. Funct. Anal.}, 268(5):1241--1276, 2015.

\bibitem{NaiboJavier2013}
Javier Duoandikoetxea and Virginia Naibo.
\newblock Mixed-norm estimates for the {$k$}-plane transform.
\newblock In {\em Excursions in harmonic analysis. {V}olume 2}, Appl. Numer.
  Harmon. Anal., pages 211--228. Birkh\"auser/Springer, New York, 2013.

\bibitem{E}
{M}.~{B}urak {E}rdo\u{g}an.
\newblock Mixed-norm estimates for a restricted {X}-ray transform in
  $\mathbb{R}^4$ and $\mathbb{R}^5$.
\newblock {\em Internat. Math. Res. Notices}, 2001(11), 2001.

\bibitem{ErdoganOberlin2010}
M.~Burak {E}rdo\u{g}an and Richard Oberlin.
\newblock Estimates for the {$X$}-ray transform restricted to 2-manifolds.
\newblock {\em Rev. Mat. Iberoam.}, 26(1):91--114, 2010.

\bibitem{Flock2016}
Taryn~C. Flock.
\newblock Uniqueness of extremizers for an endpoint inequality of the
  {$k$}-plane transform.
\newblock {\em J. Geom. Anal.}, 26(1):570--602, 2016.

\bibitem{Flockthesis}
Taryn~Cristina Flock.
\newblock {\em On extremizers for certain inequalities of the k-plane transform
  and related topics}.
\newblock ProQuest LLC, Ann Arbor, MI, 2014.
\newblock Thesis (Ph.D.)--University of California, Berkeley.

\bibitem{gel}
I.M. Gelfand, S.G. Gindikin, and M.I. Graev.
\newblock {\em Selected Topics in Integral Geometry}.
\newblock Translations of mathematical monographs. American Mathematical Soc.,
  2003.

\bibitem{GreenleafSeeger1998}
Allan Greenleaf and Andreas Seeger.
\newblock Fourier integral operators with cusp singularities.
\newblock {\em Amer. J. Math.}, 120(5):1077--1119, 1998.

\bibitem{Gressman2025}
Philip Gressman.
\newblock Testing conditions for multilinear {R}adon-{B}rascamp-{L}ieb
  inequalities.
\newblock {\em Trans. Amer. Math. Soc.}, 378(2):751--804, 2025.

\bibitem{Gressman2007}
Philip~T. Gressman.
\newblock Sharp {$L^p$}-{$L^q$} estimates for generalized {$k$}-plane
  transforms.
\newblock {\em Adv. Math.}, 214(1):344--365, 2007.

\bibitem{Gressman2009}
Philip~T. Gressman.
\newblock {$L^p$}-improving properties of averages on polynomial curves and
  related integral estimates.
\newblock {\em Math. Res. Lett.}, 16(6):971--989, 2009.

\bibitem{Gressman2013}
Philip~T. Gressman.
\newblock Uniform sublevel {R}adon-like inequalities.
\newblock {\em J. Geom. Anal.}, 23(2):611--652, 2013.

\bibitem{Gressmanintermediate2019}
Philip~T. Gressman.
\newblock Generalized curvature for certain {R}adon-like operators of
  intermediate dimension.
\newblock {\em Indiana Univ. Math. J.}, 68(1):201--246, 2019.

\bibitem{Gressman2019}
Philip~T. Gressman.
\newblock {On the Oberlin affine curvature condition}.
\newblock {\em Duke Mathematical Journal}, 168(11):2075 -- 2126, 2019.

\bibitem{Gressman2022}
Philip~T. Gressman.
\newblock Geometric averaging operators and nonconcentration inequalities.
\newblock {\em Anal. PDE}, 15(1):85--122, 2022.

\bibitem{Gressman2024}
Philip~T. Gressman.
\newblock Generalized sublevel estimates for form-valued functions and related
  results for radon-like transforms, 2024.

\bibitem{JF}
Fritz John.
\newblock The ultrahyperbolic differential equation with four independent
  variables.
\newblock {\em Duke Math. J.}, 4(2):300--322, 06 1938.

\bibitem{kattao}
Nets Katz and Terence Tao.
\newblock Recent progress on the kakeya conjecture, 2000.

\bibitem{labatao}
Izabella Laba and Terence Tao.
\newblock An {X}ray transform estimate in $\mathbb{R}^n$.
\newblock {\em Revista Matemática Iberoamericana}, 17(2):375--407, 2001.

\bibitem{NL}
{N}orberto {L}aghi.
\newblock A note on restricted {X}-ray transforms.
\newblock {\em Math. Proc. {C}ambridge Philos. Soc.}, 146(3):719--729, 05 2009.

\bibitem{LL}
E.H. Lieb and M.~Loss.
\newblock {\em Analysis}.
\newblock Crm Proceedings \& Lecture Notes. American Mathematical Society,
  2001.

\bibitem{Littman1971}
Walter Littman.
\newblock {$L\sp{p}-L\sp{q}$}-estimates for singular integral operators arising
  from hyperbolic equations.
\newblock volume XXIII of {\em Proc. Sympos. Pure Math.}, pages 479--481. Amer.
  Math. Soc., Providence, RI, 1973.

\bibitem{NagelSteinWanger1985}
Alexander Nagel, Elias~M. Stein, and Stephen Wainger.
\newblock Balls and metrics defined by vector fields {I} : Basic properties.
\newblock {\em Acta Mathematica}, 155:103 -- 147, 1985.

\bibitem{Oberlin1987}
Daniel~M. Oberlin.
\newblock Convolution estimates for some measures on curves.
\newblock {\em Proc. Amer. Math. Soc.}, 99(1):56--60, 1987.

\bibitem{Oberlin2000}
Daniel~M. Oberlin.
\newblock An estimate for a restricted {X}-ray transform.
\newblock {\em Canad. Math. Bull.}, 43(4):472--476, 2000.

\bibitem{PramanikSeeger2006}
Malabika Pramanik and Andreas Seeger.
\newblock {$L^p$} {S}obolev regularity of a restricted {X}-ray transform in
  {$\Bbb R^3$}.
\newblock In {\em Harmonic analysis and its applications}, pages 47--64.
  Yokohama Publ., Yokohama, 2006.

\bibitem{PramanikSeeger2007}
Malabika Pramanik and Andreas Seeger.
\newblock {$L^p$} regularity of averages over curves and bounds for associated
  maximal operators.
\newblock {\em Amer. J. Math.}, 129(1):61--103, 2007.

\bibitem{PramanikSeeger2021}
Malabika Pramanik and Andreas Seeger.
\newblock {$L^p$}-{S}obolev estimates for a class of integral operators with
  folding canonical relations.
\newblock {\em J. Geom. Anal.}, 31(7):6725--6765, 2021.

\bibitem{Rubin2004}
Boris Rubin.
\newblock Radon transforms on affine {G}rassmannians.
\newblock {\em Trans. Amer. Math. Soc.}, 356(12):5045--5070, 2004.

\bibitem{Rubin2019}
Boris Rubin.
\newblock Norm estimates for {$k$}-plane transforms and geometric inequalities.
\newblock {\em Adv. Math.}, 349:29--55, 2019.

\bibitem{Rubin2023}
Boris Rubin.
\newblock Fractional integrals associated with {R}adon transforms.
\newblock {\em Forum Math.}, 35(6):1727--1759, 2023.

\bibitem{Schlag1997}
Wilhelm Schlag.
\newblock A generalization of {B}ourgain's circular maximal theorem.
\newblock {\em J. Amer. Math. Soc.}, 10(1):103--122, 1997.

\bibitem{Schlag1998kakeya}
Wilhelm Schlag.
\newblock A geometric inequality with applications to the {K}akeya problem in
  three dimensions.
\newblock {\em Geom. Funct. Anal.}, 8(3):606--625, 1998.

\bibitem{Schlag1998}
Wilhelm Schlag.
\newblock A geometric proof of the circular maximal theorem.
\newblock {\em Duke Math. J.}, 93(3):505--533, 1998.

\bibitem{Seeger1998}
Andreas Seeger.
\newblock Radon transforms and finite type conditions.
\newblock {\em Journal of the American Mathematical Society}, 11(4):869 -- 897,
  1998.

\bibitem{Plamen-Uhlmann2012}
Plamen Stefanov and Gunther Uhlmann.
\newblock The geodesic {X}-ray transform with fold caustics.
\newblock {\em Anal. PDE}, 5(2):219 -- 260, 2012.

\bibitem{Stovall2009}
Betsy Stovall.
\newblock Endpoint bounds for a generalized {R}adon transform.
\newblock {\em J. Lond. Math. Soc. (2)}, 80(2):357--374, 2009.

\bibitem{Stovall2010}
Betsy Stovall.
\newblock Endpoint {$L^p\to L^q$} bounds for integration along certain
  polynomial curves.
\newblock {\em J. Funct. Anal.}, 259(12):3205--3229, 2010.

\bibitem{Stovall2011}
Betsy Stovall.
\newblock {$L^p$} improving multilinear {R}adon-like transforms.
\newblock {\em Rev. Mat. Iberoam.}, 27(3):1059--1085, 2011.

\bibitem{Stovall2020}
Betsy Stovall.
\newblock Extremizability of {F}ourier restriction to the paraboloid.
\newblock {\em Adv. Math.}, 360:106898, 18, 2020.

\bibitem{Tao-Wirght2001}
Terence Tao and Jim Wright.
\newblock ${L}^p$ improving bounds for averages along curves.
\newblock {\em Journal of the American Mathematical Society}, 16:605 -- 638,
  2001.

\bibitem{triebel}
Hans Triebel.
\newblock {\em Interpolation Theory - Function Spaces - Differential
  Operators}.
\newblock Wiley, 1999.

\bibitem{Wolff1995}
Thomas Wolff.
\newblock An improved bound for {K}akeya type maximal functions.
\newblock {\em Rev. Mat. Iberoamericana}, 11(3):651--674, 1995.

\bibitem{Wolff1997}
Thomas Wolff.
\newblock A {K}akeya-type problem for circles.
\newblock {\em Amer. J. Math.}, 119(5):985--1026, 1997.

\bibitem{Wolf01}
Thomas Wolff.
\newblock A sharp bilinear cone restriction estimate.
\newblock {\em Ann. of Math. (2)}, 153(3):661--698, 2001.

\end{thebibliography}


\end{document}